%% file: cayleyhom-final.tex
\documentclass[12pt,a4paper]{amsart}

\usepackage[T1]{fontenc}
\usepackage{a4wide}
\usepackage{amsmath,amsthm,amssymb}
\usepackage{bbm}
\usepackage{paralist}
\usepackage{color,graphicx}
\usepackage{ae,aecompl}

\newtheorem{theorem}{Theorem}[section]
\newtheorem{proposition}[theorem]{Proposition}

\theoremstyle{definition}
\newtheorem{definition}[theorem]{Definition}
\newtheorem{remark}[theorem]{Remark}
\newtheorem{observation}[theorem]{Observation}
\newtheorem{example}[theorem]{Example}

\newcommand{\CC}{\mathcal{C}}
\newcommand{\DD}{\text{\sffamily\upshape D}}
\newcommand{\FF}{\mathcal{F}}

\newcommand{\LL}{\mathcal{L}}
\newcommand{\PP}{\mathcal{P}}

\newcommand{\NN}{\mathbbm{N}}
\newcommand{\ZZ}{\mathbbm{Z}}
\newcommand{\RR}{\mathbbm{R}}
\renewcommand{\aa}{\boldsymbol{a}}
\newcommand{\xx}{\boldsymbol{x}}
\newcommand{\ba}{\textstyle\boldsymbol{\frac{1}{a}}}
\newcommand{\bg}{\textstyle\boldsymbol{\frac{1}{g}}}
\newcommand{\llambda}{\textstyle\boldsymbol{\lambda}}

\newcommand{\xmapsto}[1]{\mapstochar\xrightarrow{#1}}
\newcommand*{\downmapsto}{%
  {\mathrel
    {
      {\offinterlineskip
      \mbox{$\mathchar'26$}\kern-2mm\raisebox{-0.5mm}{\mbox{$\downarrow$}}
    }}}%
}

\DeclareMathOperator{\conv}{conv}
\DeclareMathOperator*{\join}{\star}
\DeclareMathOperator{\cross}{Cr}
\DeclareMathOperator{\Aut}{Aut}
\DeclareMathOperator{\Ind}{Ind}
\DeclareMathOperator{\IInd}{\text{\sffamily\upshape Ind}}
\DeclareMathOperator{\Cl}{\text{\sffamily\upshape Cl}}
\DeclareMathOperator{\TT}{\text{\sffamily\upshape T}}
\DeclareMathOperator{\vertices}{vert}
\DeclareMathOperator{\sk}{sk}
\newcommand{\Hom}{\text{\sffamily\upshape Hom}}
\newcommand{\IHom}{\text{\upshape{I}\sffamily\upshape Hom}}
\DeclareMathOperator{\IC}{\text{\sffamily\upshape IC}}
\DeclareMathOperator{\id}{id}
\DeclareMathOperator{\supp}{supp}

\graphicspath{{graphics/}}
\makeatletter
\def\input@path{{graphics/}}    
\makeatother

\title{Dissections, $\Hom$-complexes and the Cayley trick}

\author{Julian Pfeifle}

\address{Departament de Matem\`{a}tica Aplicada II, Universitat
  Polit\`{e}cnica de Catalunya}

\email{julian.pfeifle@upc.edu}

\thanks{The author was supported by a \emph{Juan de la Cierva} grant from the
  Spanish Ministry of Education (MEC), and by projects MCYT BFM2003-00368 and
  MTM2005-08618-C02-01.}

\date{June 27, 2006}


\begin{document} 

\begin{abstract}
  We show that certain canonical realizations of the complexes $\Hom(G,H)$ and
  $\Hom_+(G,H)$ of (partial) graph homomorphisms studied by Babson and Kozlov
  are in fact instances of the polyhedral Cayley trick. For $G$~a complete
  graph, we then characterize when a canonical projection of these complexes
  is itself again a complex, and exhibit several well-known objects that arise
  as cells or subcomplexes of such projected $\Hom$-complexes: the dissections
  of a convex polygon into $k$-gons, Postnikov's generalized permutohedra,
  staircase triangulations, the complex dual to the lower faces of a cyclic
  polytope, and the graph of weak compositions of an integer into a fixed
  number of summands.
\end{abstract}

\maketitle

\section{Introduction}

A \emph{homomorphism} from a graph $G$ to a graph~$H$ is a map
$\varphi:V(G)\to V(H)$ between their vertex sets such that
$(\varphi(x),\varphi(y))$ is an edge of~$H$ whenever $(x,y)$ is an
edge of~$G$. The \emph{$\Hom$-complex} $\Hom(G,H)$ is a polytopal
complex associated to the set of all homomorphisms from~$G$ to~$H$
that, intuitively, collects ``compatible'' homomorphisms into
polytopal cells.  

Recently, the study of $\Hom$-complexes of graphs has led to a number
of successes in topological combinatorics.  One example is the recent
proof of the \emph{Lovász Conjecture} by Babson and
Kozlov~\cite{Babson-Kozlov04} (with simplications and extensions by
Schultz~\cite{Schultz05}; see also the excellent survey
article~\cite{Kozlov05}). This result provides a lower bound for the
chromatic number of a graph~$G$ in terms of a topological property
(the connectivity) of the associated $\Hom$-complex $\Hom(C,G)$, where
$C$~is an odd cycle. In the course of the original proof of the Lovász
Conjecture, Babson and Kozlov define a certain simplicial complex
$\Hom_+(G,H)$, which is related to the set of all ``partial
homomorphisms'' from~$G$ to~$H$, i.e., homomorphisms from an induced
subgraph of~$G$ to the graph~$H$. The definitions of~$\Hom(G,H)$
and~$\Hom_+(G,H)$ are purely combinatorial.

One of the first goals of this paper is to show that certain canonical
realizations of the complexes $\Hom(G,H)$ and~$\Hom_+(G,H)$ in
Euclidean space are related via a (by now) rather famous geometric
construction, namely the \emph{polyhedral Cayley Trick} due to
Sturmfels, Huber, Rambau and
Santos~\cite{Sturmfels94}~\cite{Huber-Rambau-Santos00}~\cite{Santos05}.
This is done in Theorem~\ref{thm:intersect}.

Next, we use the canonical geometric embedding of these complexes to
project them, again canonically, to a lower-dimensional subspace. In
general, this projection $\pi\Hom(G,H)$ is not itself a polytopal
complex because the projected cells need not intersect in common
faces. However, we can characterize the shape of these projected cells
(Theorem~\ref{thm:permutohedron}): They are exactly the
\emph{generalized permutohedra} found by Postnikov~\cite{Postnikov05}.

In view of our application to dissection complexes, we then concentrate on the
special case $G=K_g$. (Note that in the literature on topological methods in
graph coloring, one is usually interested in the case $H=K_h$.) In this case,
we can characterize when the projected $\Hom$-complexes are themselves
polytopal or simplicial complexes (Theorem~\ref{thm:hypersimplex}): This
happens if and only if $\omega(H)=g$, which means that the number of vertices
in a largest clique of~$H$ is~$g$. Along the way, we define two more complexes
associated to $\Hom$-complexes, namely \emph{transversal} complexes
$\Hom^t(G,H)$ and \emph{induced} ones, $\IHom(G,H)$; moreover, we show that
for any graph~$H$, the $1$-skeleton of the projection $\pi\Hom(K_g,H)$ is a
subcomplex of the $1$-skeleton of a hypersimplex
(Proposition~\ref{prop:hypersimplex}).

We are now ready to apply these tools to \emph{dissection complexes}.  For
this, consider the set of dissections of a convex $(m(k-2)+2)$-gon into
$m$~convex $k$-gons.\footnote{As a historical aside, the interest in these
  objects goes back at least to~1791, when they were studied by Euler's
  assistant and student Nikolaus Fuss~\cite{Fuss91} in St.~Petersburg; cf.\ 
  also~\cite{Przytycki-Sikora00}} We denote by~$\delta(k,m)$ the set of all
diagonals that can arise in such a dissection, and by~$I(k,m)$ the
\emph{independence graph} on the vertex set~$\delta(k,m)$, i.e., we connect
two diagonals by an edge if the relative interiors of the diagonals do not
intersect. In Proposition~\ref{prop:hom}, we find some old acquaintances
inside the projected complexes $\DD(k,m)=\pi\Hom\big(K_{m-1},I(k,m)\big)$ and
$\DD_+(k,m)=\pi\Hom_+\big(K_{m-1},I(k,m)\big)$. Namely, the simplicial complex
induced on the set of transversal $(m-2)$-dimensional faces of~$\DD_+(k,m)$ is
a simplicial complex~$\TT(k,m)$ already considered by
Tzanaki~\cite{Tzanaki05}, and the $1$-skeleton of~$\DD(k,m)$ is the \emph{flip
  graph} on the dissections considered in~\cite{Huemer-Hurtado-Pfeifle05}.

Finally, in Section~\ref{sec:staircase} we prove interesting
isomorphisms between a certain polytopal complex~$\CC(r,s)$ whose
graph is the graph of all \emph{weak compositions} of the positive
integer~$r$ into $s$~non-negative summands, a certain induced
subcomplex of a \emph{polar-to-cyclic} polytope, and the
\emph{staircase triangulation}~$\Sigma(r,s)$ of the product of
simplices $\Delta^{r-1}\times\Delta^{s-1}$ --- of course, here the
polyhedral Cayley trick again plays a key role. As our last result, we
show that~$\CC(r,s)$ and $\Sigma(r,s)$ are basic building
blocks of~$\DD(k,m)$, respectively~$\DD_+(k,m)$.

\section{The Cayley trick and $\Hom$-complexes}

\subsection{The polyhedral Cayley trick}

Let $e_1,\dots, e_a$ be a linear basis of $\RR^a$.

\begin{proposition}\label{prop:cayley}
  \cite{Sturmfels94}~\cite{Huber-Rambau-Santos00}~\cite{Santos05} Fix
  real $\lambda_1,\dots,\lambda_n$ such that each $\lambda_i>0$ and
  $\sum_{i=1}^n\lambda_i=1$. Then, for any polytopes $P_1, \dots,
  P_n\subset\RR^d$, there is an isomorphism between the posets of
  polyhedral subdivisions of the \emph{Cayley embedding}
  \[ 
     \CC(P_1,\dots, P_n)  \ = \  \conv\ \bigcup_{i=1}^n P_i\times e_i 
     \ \subset \ \RR^d\times \RR^n
  \]
  of the $P_i$'s and the poset of mixed subdivisions of the Minkowski sum
  $\sum_{i=1}^n \lambda_i P_i\subset\RR^d$, both ordered by refinement. The
  bijection between two corresponding subdivisions is given by intersecting a
  polyhedral subdivision $\mathcal{P}$ of $\CC(P_1,\dots,P_n)$ with the
  $d$-dimensional plane $L=\RR^d\times(\lambda_1,\dots, \lambda_n)$. This
  intersection produces from each cell $\conv \bigcup_{i=1}^n Q_i\times e_i$
  of~$\mathcal{P}$ the weighted Minkowski sum $(\sum_{i=1}^n \lambda_i
  Q_i)\times(\lambda_1,\dots,\lambda_n)$, where $Q_i\subset P_i$ are
  subpolytopes.  \hfill$\Box$
\end{proposition}

\subsection{Simultaneous instances of the Cayley trick, related by
  joins and projections}

To paraphrase Proposition~\ref{prop:cayley}, the Cayley trick relates a
``Cayley object''~---~namely a cell~$Q=\conv \bigcup_{i=1}^n Q_i\times e_i$ of
a polyhedral subdivision~$\PP$ of the Cayley embedding of the polytopes
$P_1$,~\dots,~$P_n$~---~to its corresponding ``Minkowski object'', namely the
(weighted and embedded) Minkowski sum of the subpolytopes $Q_1$,~\dots,~$Q_n$
of the~$P_i$'s. The agent that produces this correspondence is a ``morphing
plane''~$L$ that intersects the cells of the subdivision (and also determines
the weights in the Minkowski sum; however, here we will mostly just need the
case of equal weights $\lambda_1=\dots=\lambda_n=\frac{1}{n}$).

In this paper, we will in fact work with \emph{two} simultaneous
``horizontal'' instances of the Cayley trick, which will be related to each
other by a ``vertical'' projection called~$\pi_\Box$. The ``bottom'' instance
of the Cayley trick will be much as we have just outlined, but the ``top''
instance will be rather special: The top Cayley objects will always be joins
of simplices (labeled by~``$J$''), and the top Minkowski objects will be
products of polytopes (labeled by~``$\Pi$''); similarly, we label the bottom
Cayley objects by~``$C$'' and the bottom Minkowski objects
by~``$M$''.\footnote{Thanks to one of the referees for suggesting this
  language.}  We summarize this situation in the commutative diagram
\[
\begin{matrix}
  {}^J\! Q & \xmapsto{\ \iota_L\ } &
  {}^\Pi\!\big((\sum_{i=1}^n\lambda_i Q_i) 
  \times\llambda\big) \\[2ex]
  \quad\;\;\downmapsto\pi_\Box && \downmapsto\pi_\Box\\[2ex]  
  {}^C\!\!\left(\pi_\Box Q\right) 
  & \xmapsto{\ \iota_{\pi_\Box(L)}\ }
  &   {}^M\!\big((\sum_{i=1}^n\lambda_i \pi_\Box(Q_i))
  \times\llambda\big) ,
\end{matrix}
\]
where for any union of polytopes $\PP\subset\RR^d$ and any affine subspace
$K\subset\RR^d$, we denote the intersection of~$\PP$ with~$K$ by
$\iota_K(\PP)=\PP\cap K$, and $\llambda=(\lambda_1,\dots,\lambda_n)$. The
$\Hom$-complexes central to this paper, and their various projections, fit
roughly as follows into this diagram:
\[
\begin{matrix}
  {}^J\!\!\left(\Hom_+(G,H)\right) & \xmapsto{\ \iota_L\ } &
  {}^\Pi\!\left(\Hom(G,H)\times\llambda\right) \\[2ex]
  \quad\;\;\downmapsto\pi_\Box && \downmapsto\pi_\Box\\[2ex]  
  {}^C\!\big(\pi_\Box\Hom_+(G,H)\big) 
  & \xmapsto{\ \iota_{\pi_\Box(L)}\ }
  &   {}^M\!\big(\pi_\Box\Hom(G,H)\times\llambda\big)  \\[2ex]
   \quad\;\;\downmapsto\pi_\Delta && \downmapsto\pi_\Delta\\[2ex]
   \DD_+(k,m) && \DD(k,m)\,.
\end{matrix}
\]
To aid the intuition of the reader, we have also included the dissection
complexes associated to the special case $G=K_{m-1}$, $H=I(k,m)$ into this
sketch; the projections $\pi_\Box$ and $\pi_\Delta$ will be defined in a
minute.  

\subsubsection{The top instance}

To explicitly define the objects participating in the ``top'' instance
of the Cayley trick, we first assemble some notation. For sets $A$ and
$B$ of respective cardinalities $a=|A|$ and $b=|B|$, denote by
$\Delta_A$ the simplex $\conv\{e_i:i\in A\}\subset\RR^{|A|}$ on the
vertex set~$A$, so that $\dim \Delta_A=|A|-1$, and similarly
for~$\Delta_B$. We will often not distinguish between a subset
$\tau\subset B$ and a face $\tau$ of~$\Delta_B$.

A key observation is now that the abstract $(ab-1)$-dimensional
simplex that arises as  the join $\join_{x\in A}\Delta_{B}$ can be
geometrically realized as the Cayley embedding of the polytopes
$\mu_1(\Delta_B)$, \dots, $\mu_a(\Delta_B)$ into
$\RR^{ab}\times\RR^a$, where $\mu_i:\RR^b\hookrightarrow\RR^{ab}$ is
the inclusion of $\RR^b$ into the $i$-th component of
$\RR^{ab}=\RR^b\times\dots\times\RR^b$.  We obtain
\[
   \join_{i\in A}\Delta_{B} \ = \ 
   \CC\big(\mu_1(\Delta_B),\dots,\mu_a(\Delta_B)\big) \ = \
   \conv\ \bigcup_{i=1}^a\mu_i(\Delta_{B})\times e_i\,.
\]
Observe that this Cayley embedding is indeed a simplex (and therefore equal to
the join $\join_{x\in A}\Delta_B$), because the $\mu_i(\Delta_B)$ are affinely
independent from each other. Moreover, all faces of~$\join_{i\in A}\Delta_B$
are of the form
\[
   \sigma\ = \ \join_{i\in  A}\sigma_i
   \ = \ \CC\big(\mu_1(\sigma_1),\dots,\mu_a(\sigma_a)\big)
\]
for some collection of faces $(\sigma_i:i\in A)$ of $\Delta_B$. In accordance
with our earlier discussion, we will sometimes explicitly identify such a face
$\sigma={}^J\!\sigma$ as being of ``join type''.

Similarly, the Minkowski object
$\frac1a(\mu_1(\sigma_1)+\dots+\mu_a(\sigma_a))$ corresponding to
${}^J\!\sigma$ is in fact a cartesian product
$\frac1a(\mu_1(\sigma_1)\times\dots\times\mu_a(\sigma_a))$, because
the $\mu_i(\sigma_i)$ lie in mutually skew subspaces by construction;
hence we will refer to this Minkowski object as being of ``product
type''~$\Pi$.

\subsubsection{The projections}

Next, we define the two projections
\begin{eqnarray*}
\pi_\Box \ :\ \RR^{ab}\times\RR^a &\to& \RR^b\times\RR^a,\\
\pi_\Delta \ : \ \RR^{b\phantom{a}}\times\RR^a &\to& \RR^b, 
\end{eqnarray*}
as follows. The map $\pi_\Delta$ is just
the projection onto the first factor; its purpose is to eliminate the
extraneous factor ``$\times\llambda$''.  As for $\pi_\Box$, on the one hand we
want it to leave the last factor~$\RR^a$ (and in particular the
point~$\llambda$) invariant; on the other, for reasons that will become clear
below, we would like it to superimpose all $a$~copies of~$\RR^b$ in the
factor~$\RR^{ab}$ onto each other. Therefore, we choose the matrix
of~$\pi_\Box$ to be
\[
    \begin{pmatrix}
      \mathbbm{1}_b & \cdots  & \mathbbm{1}_b & 0\\
      0 & \cdots & 0 &  \mathbbm{1}_a
    \end{pmatrix},
\]
where the $\mathbbm{1}_k$ denote $k\times k$ unit matrices, and the
zeros stand for null matrices of the appropriate size.  Note that,
loosely speaking, each $\mu_i$~is a section of~$\pi_\Box$, in the
sense that $\pi_\Box|_{\RR^{ab}}\circ\mu_i=\id_{\RR^b}$ for
all~$i$. In particular,
\begin{equation}\label{eq:projcayley}
  \pi_\Box(\sigma) \ = \
   \pi_\Box\;\CC\big(\mu_1(\sigma_1), \dots, \mu_a(\sigma_a)\big)
   \ = \ \CC(\sigma_1,\dots,\sigma_a)
\end{equation}
for any face $\sigma=\star_{i\in A}\sigma_i$ of $\join_{i\in A}\Delta_B$.

Finally, let us fix our ``morphing plane''~$L$ once and for all as the
$ab$-dimensional plane $L\subset\RR^{ab}\times\RR^a$ defined by
\begin{equation}\label{eq:L}
   L \ =  \   \RR^{ab}\times\ba\,,
\end{equation}
where here and throughout we set
$\ba=(\frac{1}{a},\dots,\frac{1}{a})\in\RR^a$, so that
$\pi_\Box(L)=\RR^b\times\ba$. 

We summarize our discussion in the following proposition.

\begin{proposition}\label{prop:diagram}
  Let $\sigma = {}^J\!\sigma = \star_{i\in
    A}\sigma_i=\CC\big(\mu_1(\sigma_1),\dots,\mu_a(\sigma_a)\big)$ be
  a face of $\join_{i\in A}\Delta_B$.  Then the following diagram
  commutes:
  \smallskip
  \[
  \begin{matrix}
    \displaystyle\join_{i\in A}\Delta_B 
    &\supset& 
       \sigma_1 \star \dots \star \sigma_a
    & \xmapsto{\ \iota_L \ } 
    & \tfrac{1}{a}\big(\mu_1(\sigma_1)\times \dots \times \mu_a(\sigma_a)\big) 
      \times\ba
    &\subset& \Delta_{A\times B}\times\ba  \\[2ex]
    &&\downmapsto\pi_\Box && \downmapsto\pi_\Box\\[2ex]
    \Delta_B\times\Delta_A
    &\supset&
    \CC(\sigma_1,\dots,\sigma_a\big)
    & \xmapsto{\ \iota_{\pi_\Box(L)} \ } 
    & \tfrac{1}{a}\big(\sigma_1+ \dots + \sigma_a\big) \times\ba 
    & \subset & \Delta_B\times\ba\\[2ex]
    &&\downmapsto\pi_\Delta && \downmapsto\pi_\Delta\\[2ex]  
    \Delta_B &\supset&
    \conv\ \bigcup_{i=1}^a \sigma_i && 
    \tfrac{1}{a}\big(\sigma_1+ \dots + \sigma_a\big) 
    & \subset & \Delta_B
  \end{matrix}
  \]
  The (reverse) inclusions on the left-hand side of the
  diagram map vertices to vertices. This is generally not the case for the
  inclusions on the right-hand side.

\end{proposition}

\begin{proof}
  It suffices to check the top middle square of the diagram. The
  horizontal maps are well-defined because they are just applications
  of the polyhedral Cayley trick, and the vertical maps are
  well-defined by~\eqref{eq:projcayley} and the linearity
  of~$\pi_\Box$. Taken together, this also proves commutativity.  The
  rest of the diagram follows by checking the definitions. 
\end{proof}

\begin{observation}
  If $\Delta_B$ is a join $\Delta_B=\join_{j\in C}\Delta_D$ and each
  $\sigma_i$ resides in a different copy of~$\Delta_D$ (which in particular
  implies $\sigma_i\cap\sigma_j=\emptyset$ for all $i\ne j$), we can glue the
  first row of another copy of this diagram onto the last row of this one, and
  in particular fill in the missing map in the last row.
\end{observation}

\begin{proof}
  Include $\Delta_B\subset\RR^b=\RR^{cd}$ into $\RR^{cd}\times\RR^c$, where
  $c=|C|$ and $d=|D|$, by the map that sends the $j$-th block of variables,
  $(x_{(j-1)d+1},\dots,x_{jd})$, of~$\RR^{cd}$ to the block
  $(x_{(j-1)d+1},\dots,x_{jd}$, $1-x_{(j-1)d+1}-\dots-x_{jd})$
  of~$\RR^{cd}\times\RR^c$, for $1\le j\le c$. This brings $\join_{j\in
    C}\Delta_D$ into the required canonical form.
\end{proof}

See Example~\ref{ex:detail} below for a detailed calculation with coordinates;
here we first present a more conceptual illustration.

\begin{figure}[htbp]
  \centering
  \input{graphics/cayley1.pstex_t}
  \caption{Two instances of the projection~$\pi_\Delta$ from
    Proposition~\ref{prop:diagram}, applied to faces
    $\CC(\sigma_1,\sigma_2)\subset\Delta_{\{3,4,5\}} \times\Delta_{\{1,2\}}$.
    \emph{Left:} $\sigma_1=\{3,4\}$ and $\sigma_2=\{4,5\}$, so that
    $\sigma_1\cap\sigma_2\ne\emptyset$. \emph{Right:} $\sigma_1=\{3,4\}$ and
    $\sigma_2=\{5\}$, so that $\sigma_1\cap\sigma_2=\emptyset$. In this case,
    $\Delta_B=\Delta_{\{3,4,5\}}=\sigma_1\star\sigma_2$, and we obtain
    $\pi_\Delta\pi_\Box\,\iota_L(\sigma) =
    \frac12(\sigma_1+\sigma_2)\subset\Delta_B$ via the polyhedral Cayley
    trick, by intersecting $\pi_\Delta\pi_\Box(\sigma) =
    \conv(\sigma_1\cup\sigma_2)\subset\Delta_B$ with the affine subspace~$L'$
    on the bottom right.}
  \label{fig:cayley1}
\end{figure}
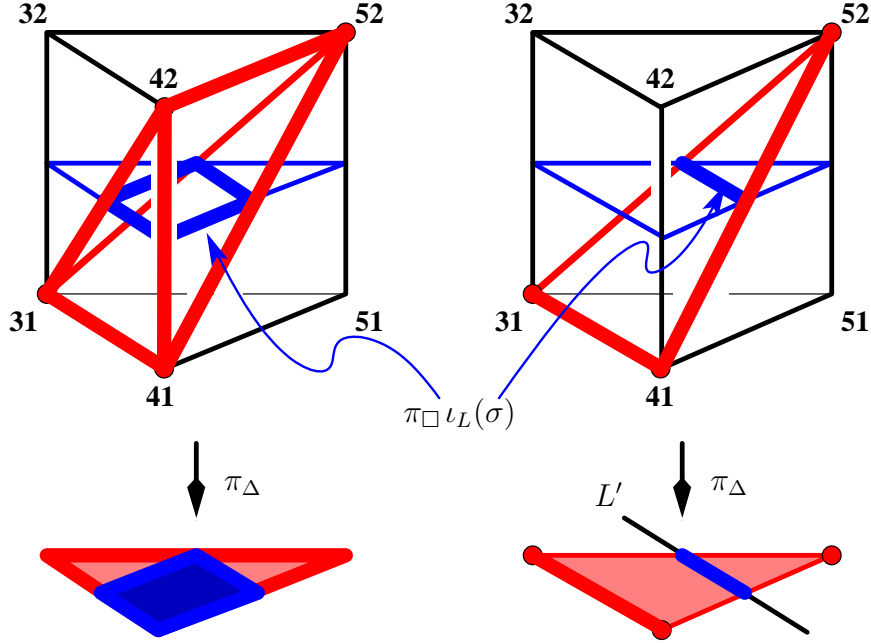

\begin{example}
  Let $A=\{1,2\}$, $B=\{3,4,5\}$, $\sigma_1=\{3,4\}$ and $\sigma_2=\{4,5\}$,
  and let us evaluate the middle row (i.e., the ``lower instance of the Cayley
  trick'') of the preceding diagram. We see that $\Delta_B\times\Delta_A$ is a
  triangular prism with vertex set $B\times A$ (but embedded into
  $\RR^3\times\RR^2$), the Cayley embedding $\CC(\sigma_1,\sigma_2)$ is the
  tetrahedron $T=\conv\{31,41,42,52\}$, and the corresponding cell of the
  subdivision is the quadrilateral that results from slicing~$T$ with
  $\pi_\Box(L)=\RR^3\times(\frac12,\frac12)$. When we apply~$\pi_\Delta$, on
  the bottom row of the diagram we obtain on the left-hand side
  $\conv(\sigma_1\cup\sigma_2)=\Delta_B$, and on the right-hand side
  $\frac12(\sigma_1+\sigma_2)$, a quadrilateral in~$\Delta_B$ that is the
  Minkowski sum of two edges scaled by~$\frac12$; see the left-hand side of
  Figure~\ref{fig:cayley1}.
  
  On the other hand, if we choose $\sigma_1=\{3,4\}$ and
  $\sigma_2=\{5\}$ to be disjoint, we obtain on the left-hand side
  of the bottom row again $\conv\{3,4,5\}=\Delta_B$, but on the
  right-hand side the segment $s=\frac12(34+5)$. Observe that in this
  situation, $\Delta_B=\sigma_1\star\sigma_2$, and that $s$~arises by
  applying the Cayley trick to this join. See the right-hand side of
  Figure~\ref{fig:cayley1}.
\end{example}

\subsection{$\Hom$-complexes}
Let $G$ and $H$ be graphs on $g=|V(G)|$ and $h=|V(H)|$ vertices.  When
convenient, we will identify $V(G)$ and $V(H)$ with $[g]=\{1,2,\dots,g\}$,
respectively~$[h]$.
A \emph{homomorphism} from $G$ to $H$ is a map $\varphi:V(G)\to V(H)$ such
that for any edge $(x,y)$ of~$G$, $(\varphi(x),\varphi(y))$ is an edge of~$H$.

Recall from~\cite{Kozlov05} the following material related to the set
of all homomorphisms between $G$~and~$H$.  
Let $\times_{x\in V(G)}\Delta_{V(H)}$ denote the cartesian product of
$g$~copies of the simplex~$\Delta_{V(H)}$, so that the copies
of~$\Delta_{V(H)}$ are labeled by the vertices of~$G$. Similarly,
$\join_{x\in V(G)}\Delta_{V(H)}$ is the join of $g$ labeled copies of
$\Delta_{V(H)}$.  Note that $\times_{x\in V(G)}\Delta_{V(H)}$ is a
$g(h-1)$-dimensional polytope that is a product of simplices, while
$\join_{x\in V(G)}\Delta_{V(H)}$ is a simplex of dimension~$gh-1$. We
will always think of $\join_{x\in V(G)}\Delta_{V(H)}$ as being embedded
in~$\RR^{gh}\times\RR^g$ as in Proposition~\ref{prop:diagram}.

\smallskip
The following two complexes have proved to be useful in topological
combinatorics; see~\cite{Kozlov05} for a survey.

\begin{definition} \cite{Kozlov05} 
  \begin{compactenum}[\upshape(a)] \samepage
  \item $\Hom(G,H)$ is the polytopal subcomplex of $\times_{x\in
      V(G)}\Delta_{V(H)}$ of all cells $\times_{x\in V(G)}\sigma_x$ such that
    if $(x,y)\in E(G)$, then $(\sigma_x,\sigma_y)$ is a complete bipartite
    subgraph of~$H$.
        
  \item $\Hom_+(G,H)$ is the simplicial subcomplex of $\join_{x\in
      V(G)}\Delta_{V(H)}$ of all simplices $\join_{x\in V(G)}\sigma_x$ such
    that if $(x,y)\in E(G)$ and both $\sigma_x$ and $\sigma_y$ are nonempty,
    then $(\sigma_x,\sigma_y)$ is a complete bipartite subgraph of~$H$.
  \end{compactenum}  
\end{definition}

Note that the bipartite subgraphs are not required to be induced.  Moreover,
by definition all cells of $\Hom(G,H)$ are products of simplices. We will
sometimes identify faces $\sigma={}^\Pi\!\sigma=\times_{x\in V(G)}\sigma_x$
of~$\Hom(G,H)$, respectively faces $\sigma={}^J\!\sigma = \join_{x\in
  V(G)}\sigma_x$ of~$\Hom_+(G,H)$, with the ordered list of (non-empty) labels
$(\lambda_1,\dots,\lambda_g)$, where $V(G)=[g]$ and $\lambda_i\subset V(H)$~is
the vertex set of the simplex~$\sigma_i$. Moreover, define
$L=\RR^{gh}\times\bg$ as in~\eqref{eq:L}.

\begin{definition}
  A face $\join_{x\in V(G)}\sigma_x$ is \emph{transversal} if $|\sigma_x|>0$
  for all $x\in V(G)$.  The simplicial complex $\Hom_+^t(G,H)$ is the
  subcomplex of $\Hom_+(G,H)$ induced by the set of all transversal faces.
\end{definition}

\begin{theorem}\label{thm:intersect}
  \begin{asparaenum}[\upshape(i)]\samepage
  \item $\iota_L\Hom_+^t(G,H) = \iota_L\Hom_+(G,H) = \Hom(G,H)\times\bg$. In
    particular, we obtain a canonical embedding of all these complexes
    into the same Euclidean space.
    
  \item The following diagram commutes: 
  \[
    \begin{matrix}
    \displaystyle\join_{x\in V(G)}\Delta_{V(H)} 
    &\supset& \Hom_+(G,H)
    & \xmapsto{\ \iota_L \ } 
    & \Hom(G,H) \times\bg
    &\subset& \Delta_{V(G)\times V(H)}\times\bg  \\[1.5ex]
    &&\downmapsto\pi_\Box && \downmapsto\pi_\Box\\[1.5ex]  
    \Delta_{V(H)}\times\Delta_{V(G)}
    &\supset&
    \pi_\Box\Hom_+(G,H)
    & \xmapsto{\ \iota_{\pi_\Box(L)} \ } 
    & \pi_\Box\Hom(G,H) \times\bg 
    & \subset & \Delta_{V(H)}\times\bg 
    \end{matrix}
  \]
  In particular, the image $\pi_\Box(\sigma)$ of any face~$\sigma$
  of~$\Hom_+(G,H)$ is the convex hull of some vertices
  of~$\Delta_{V(H)}\times\Delta_{V(G)}$, and $\pi_\Box\Hom(G,H) =
  \iota_{\pi_\Box(L)}\pi_\Box\Hom_+(G,H)$.

  \item The same statements hold with $\Hom_+(G,H)$ replaced by
    $\Hom_+^t(G,H)$.

  \end{asparaenum}
\end{theorem}

\smallskip
\begin{proof}
  (i) For the first equality, let $\sigma={}^J\!\sigma =\join_{x\in V(G)}\sigma_x$ be a
  simplex of~$\Hom_+(G,H)$. Since $\join_{x\in G}\Delta_{V(H)}$ is
  embedded in~$\RR^{gh}\times\RR^g$, any point $z\in\sigma$ can be
  written as the convex combination
  \begin{equation}\label{eq:sum1}
     z \ = \ \sum_{i=1}^g \lambda_i \sum_{v\in\vertices\sigma_i} \lambda'_{iv}
     \mu_i(v)\times e_i \ \in \ \RR^{gh}\times\RR^{g},
  \end{equation}
  where $\sum_{i=1}^g\lambda_i=1$ and
  $\sum_{v\in\vertices\sigma_i}\lambda'_{iv}=1$ for all $i=1,2,\dots, g$.  Now
  look at the $i$-th entry of the $\RR^g$-component of~$z$. If $|\sigma_i|>0$,
  it is $\lambda_i\sum_{v\in\vertices\sigma_i}\lambda'_{iv}=\lambda_i$, but
  otherwise it vanishes because the inner sum in~\eqref{eq:sum1} is empty.
  Therefore, no simplex of $\Hom_+(G,H)\setminus\Hom_+^t(G,H)$ intersects~$L$.
  (However, notice that not all faces of $\Hom_+^t(G,H)$ intersect~$L$.)
  
  The second equality follows from the top row of the diagram in
  Proposition~\ref{prop:diagram} with $a=g$ and $b=h$: 
 Observe that $\iota_L({}^J\!\sigma)
  = {}^\Pi(\mu_1(\sigma_1)+\dots+\mu_g(\sigma_g)) =
  {}^\Pi(\mu_1(\sigma_1)\times\dots\times\mu_g(\sigma_g))$ for any simplices
  $\sigma_1,\dots,\sigma_g\subset\Delta_{V(H)}$, because the
  $\mu_i(\sigma_i)$ lie in skew subspaces by construction.
  
  \smallskip (ii) Let $\sigma=\star_{i\in V(G)}\sigma_i$ be a face
  of~$\Hom_+(G,H)$, where $\sigma_i\subset\Delta_{V(H)}$ for all $i\in V(G)$.
  The image under~$\pi_\Box$ of
  ${}^J\!\sigma=\conv\bigcup_{i=1}^g\mu_i(\sigma_i)\times e_i$ is
  ${}^C\!(\pi_\Box(\sigma))=\conv\bigcup_{i=1}^g\sigma_i\times e_i$, which is
  the convex hull of some vertices of~$\Delta_{V(H)}\times\Delta_{V(G)}$.  The
  well-definedness and commutativity of the diagram now follows from the fact
  that $\iota_L\Hom_+(G,H) = \Hom(G,H)$ and Proposition~\ref{prop:diagram},
  after passing to faces.

  Finally, (iii) follows from (i).
\end{proof}

\subsection{Projections of $\Hom$-complexes and generalized permutohedra}

Suppose that $\Hom(G,H)$ is embedded into $\RR^{gh}\times\RR^g$ as in our
previous discussion, and denote the symmetry groups of the graphs~$G$~and~$H$
by~$S_G=\Aut(G)$, respectively $S_H=\Aut(H)$.  By~\cite{Kozlov05}, the product
group $S_G\times S_H$ acts on $\Hom(G,H)$.

Here and in the following we will use the notation
$\pi=\pi_\Delta\pi_\Box:\RR^{gh}\times\RR^g\to\RR^h$.

\begin{proposition}\label{prop:quotient} 
  $\pi\big(\gamma(\sigma)\big) = \pi(\sigma)$ for any $\gamma\in S_G$ and any
  cell~$\sigma$ of~$\Hom(G,H)$.
\end{proposition}

\begin{proof}
  Let $\sigma={}^\Pi\!\sigma=(\times_{i\in V(G)}\sigma_i)\times\bg$ be a cell
  of  $\Hom(G,H)\times\bg$ of ``product~type'', and let
  $\sigma'={}^\Pi\!(\sigma') = \big(\!\times_{i\in
    V(G)}\sigma_{\gamma(i)}\big)\times\bg$ be its image under~$\gamma\in S_G$.
  We will chase~$\sigma'$ around the diagram of Theorem~\ref{thm:intersect} to
  verify that $\pi_\Box(\sigma)=\pi_\Box(\sigma')$ in~$\pi_\Box\Hom(G,H)$;
  since $\pi_\Delta$ restricted to $\pi_\Box\Hom(G,H)$~is an isomorphism, this
  implies that their images in $\pi\Hom(G,H)$ coincide.
  
  First, let $\tilde\sigma$ and $\tilde\sigma'$~be the unique simplicial faces
  of~$\Hom_+(G,H)$ such that $\iota_L(\tilde\sigma)=\sigma$ and
  $\iota_L(\tilde\sigma')=\sigma'$.  We can write any point
  $z\in\pi_\Box(\tilde\sigma)$ 
  as a convex  combination
  \[
     z \ = \ \sum_{i\in V(G)} \lambda_i \ x_i \times e_i
  \]
  of points $x_i=\sum_{v\in\vertices\sigma_i}\lambda_{iv}v\in\sigma_i$, so
  that
  \[
     z' \ = \ \sum_{i\in V(G)} \lambda_i \ x_{\gamma(i)} \times e_i
  \]
  is the corresponding point in $\pi_\Box(\tilde\sigma')$ under the action
  of~$\gamma$. Now note that the $h$-plane $\pi_\Box(L)=\RR^h\times\bg$ only
  intersects those cells $\pi_\Box(\tilde\sigma)$ with all $\sigma_i$
  non-empty; the reason (as in the proof of Theorem~\ref{thm:intersect}) is
  that $\sigma_j=\emptyset$ forces the $j$-th component of the $\RR^g$-part of
  $z\in\RR^h\times\RR^g$ to be~zero.  Therefore, intersecting
  $\sigma$~and~$\tilde\sigma$ with~$\pi_\Box(L)$ forces
  $\lambda_i=\frac{1}{g}$ for all $1\le i\le g$, so that the images of
  $z$~and~$z'$ under the map $\iota_{\pi_\Box(L)}\circ\pi_\Box$ agree, which
  is what we wanted to show.
\end{proof}

\begin{observation}\label{obs:not_equivalence}
  If $\Aut(G)\subsetneq S_{|G|}$ is not the full symmetric group on $|G|$
  letters, then $\pi(\sigma)=\pi(\tau)$ may hold for faces $\sigma,\tau$
  of~$\Hom(G,H)$, even though $\tau$ is not of the form $\gamma(\sigma)$ for
  any $\gamma\in\Aut(G)$. Because of this, it would not be correct to say that
  each cell of $\pi\Hom(G,H)$ represents an $S_G$-equivalence class of faces
  of the polytopal complex~$\Hom(G,H)$.
\end{observation}

\begin{definition}
  $\pi\Hom(G,H)=\{\pi(\sigma):\sigma\text{ is a cell of }\Hom(G,H)\}$.
\end{definition}

This is not in general a polytopal complex, because cells need not
intersect in common faces. However, the next theorem identifies the
faces of~$\pi\Hom(G,H)$ as ``generalized permutohedra'', introduced by
Postnikov~\cite{Postnikov05}.

For this, let $\Gamma\subset K_{m,n}$ be a bipartite graph with $m$~``left''
and $n$~``right'' vertices. We agree to denote its edges by $(i,j)$, where $i$
is in the ``left'' part and $j$ in the ``right'' part of $\Gamma$.

For any such $\Gamma$, Postnikov defined the \emph{generalized permutohedron}
$P_\Gamma(\lambda_1,\dots,\lambda_m)$ to be the weighted Minkowski sum of
simplices $P_\Gamma=\lambda_1\Delta_{I_1}+\dots+\lambda_m\Delta_{I_m}$, where
$\lambda_i>0$ and $I_i=I_i(\Gamma)=\{j:(i,j)\in \Gamma\}\subset[n]$ for
$i=1,2,\dots,m$. These may be obtained via the polyhedral Cayley Trick from
the \emph{root polytopes}
$Q_\Gamma=\CC(\Delta_{I_1},\dots,\Delta_{I_m})\subset
\Delta_{[n]}\times\Delta_{[m]}$, by intersecting with the subspace
$\RR^n\times(\frac{\lambda_1}{\lambda},\dots,\frac{\lambda_m}{\lambda})$,
where $\lambda=\sum_{i=1}^m\lambda_i$ \cite[Corollary~14.6]{Postnikov05}.

\begin{theorem}\label{thm:permutohedron}
  Let $G$ and $H$ be graphs.  Then any cell of $\pi\Hom(G,H)$ is a
  generalized permutohedron, and any generalized permutohedron occurs
  as a cell of some~$\pi\Hom(G,H)$. Moreover, a cell~$\pi(\sigma)$
  of~$\pi\Hom(G,H)$ is a product of simplices if and only if
  $\sigma_i\cap\sigma_j=\emptyset$ for all $i\ne j$.
\end{theorem}

\begin{proof}
  We start with a $\Hom$-complex $\Hom(G,H)$, where $G$~and~$H$ are graphs
  with~$g$, respectively $h$~vertices. Let
  $\sigma=\sigma_1\times\dots\times\sigma_g$ be a cell of $\Hom(G,H)$,
  embedded in~$\Delta_{V(H)}$ as in Theorem~\ref{thm:intersect}. (We will
  temporarily forget about the extra factor~``$\times\bg$'' here.) Moreover,
  let $\tilde\sigma$~be the unique simplicial cell of~$\Hom_+(G,H)$
  such that $\sigma=\iota_L(\tilde\sigma)$, and denote
  by~$W=\bigcup_{i=1}^g\sigma_i\subset V(H)$ the union of all~$\sigma_i$,
  regarded as subsets of~$V(H)$.  Now define the bipartite graph
  $\Gamma\subset K_{g,|W|}$ by connecting a left vertex~$i$ to a right
  vertex~$j$ whenever $j\in\sigma_i$. This graph defines a root polytope
  $Q_\Gamma=Q_\Gamma(1,\dots,1)=\CC(\sigma_1, \dots,\sigma_g)$, and we only
  have to check that $Q_\Gamma=\pi_\Box\tilde\sigma$.  But this is clear by
  the diagram of Proposition~\ref{prop:diagram}.
  
  In the other direction, let us first suppose that we are given a generalized
  permutohedron $P_\Gamma(1,\dots,1)$ for some bipartite graph~$\Gamma\subset
  K_{m,n}$. Define $\sigma_i=\{j:(i,j)\in \Gamma\}\subset[n]$ for
  $i=1,2\dots,m$, and set $G_\Gamma=\Ind(\sigma_1,\dots,\sigma_m)$, the
  \emph{independence graph} that has the $\sigma_i$'s as vertices, and in
  which $\sigma_i$~is joined to~$\sigma_j$ by an edge precisely
  if~$\sigma_i\cap\sigma_j=\emptyset$.  Moreover, define a graph~$H_\Gamma$ on
  the vertex set~$[n]$ by adding the edges of a complete bipartite graph
  $(\sigma_i,\sigma_j)$ whenever $\sigma_i$~and~$\sigma_j$ form an edge in
  $G_\Gamma$. Checking the definitions yields that
  $\pi\Hom(G_\Gamma,H_\Gamma)$~contains a cell of the
  form~$P_\Gamma(1,\dots,1)$.  The more general case of permutohedra of the
  form $P_\Gamma(\lambda_1,\dots,\lambda_m)$ follows by constructing
  $\Hom(G_\Gamma,H_\Gamma)$ as the slice
  $\Hom_+(G_\Gamma,H_\Gamma)\cap\big(\RR^{|G_\Gamma|\cdot|H_\Gamma|} \times
  (\frac{\lambda_1}{\lambda},\dots,\frac{\lambda_m}{\lambda})\big)$, where
  $\lambda=\sum_{i=1}^m\lambda_i$.
\end{proof}

\begin{example}\label{ex:detail}
  Let $\Gamma\subset K_{3,3}$ be the bipartite graph with edge set
  $\{\bar1 1,\bar1 2, \bar2 1, \bar2 3,\bar3 2, \bar3 3 \}$, where we
  write left vertices with bars. The proof of the preceding theorem
  yields graphs $G_\Gamma$, $H_\Gamma$ such that
  $\pi\Hom(G_\Gamma,H_\Gamma)$ contains a hexagon
  $P_\Gamma(1,1,1)$.  Namely, the proof calls for setting
  $G_\Gamma=H_\Gamma=E_3$, the graph with $3$~vertices and no edges,
  so that $\Hom(G_\Gamma, H_\Gamma)=\Delta_{[3]}\times\Delta_{[3]}$
  and $S_{G_\Gamma}=S_3$, the symmetric group on $3$~letters.
  Moreover, $\sigma_{\bar1}=\{1,2\}$, $\sigma_{\bar2}=\{1,3\}$ and
  $\sigma_{\bar3}=\{2,3\}$, so that the vertices of the simplicial cell
  $\sigma=\conv\bigcup_{i=1}^3\mu_i(\sigma_{\bar i})\times e_i$ in
  $\Hom_+(G_\Gamma,H_\Gamma)$ are the rows of the following matrix:
  \[\small
  \begin{tabular}[c]{*{15}{c}}
    1 & 0 & 0 && 0 & 0 & 0 && 0 & 0 & 0 && 1 & 0 & 0\\
    0 & 1 & 0 && 0 & 0 & 0 && 0 & 0 & 0 && 1 & 0 & 0\\
    0 & 0 & 0 && 1 & 0 & 0 && 0 & 0 & 0 && 0 & 1 & 0\\
    0 & 0 & 0 && 0 & 0 & 1 && 0 & 0 & 0 && 0 & 1 & 0\\
    0 & 0 & 0 && 0 & 0 & 0 && 0 & 1 & 0 && 0 & 0 & 1\\
    0 & 0 & 0 && 0 & 0 & 0 && 0 & 0 & 1 && 0 & 0 & 1
  \end{tabular}
  \]
  Intersecting the convex hull of these points (a $5$-dimensional simplex
  in~$\RR^{12}$) with the $9$-dimensional plane
  $L=\RR^{3\cdot3}\times(\frac13,\frac13,\frac13)$ yields the $3$-cube cell
  $\sigma\cap L$ of~$\Hom(G_\Gamma,H_\Gamma)\times\boldsymbol{\frac13}$ whose
  vertices are $\frac13$~times the row vectors
  \[
  \small
  \begin{tabular}[c]{*{15}{c}}
    1 & 0 & 0 && 1 & 0 & 0 && 0 & 0 & 1 && 1 & 1 & 1\\ 
    1 & 0 & 0 && 0 & 0 & 1 && 0 & 0 & 1 && 1 & 1 & 1\\ 
    0 & 1 & 0 && 1 & 0 & 0 && 0 & 0 & 1 && 1 & 1 & 1\\ 
    0 & 1 & 0 && 0 & 0 & 1 && 0 & 0 & 1 && 1 & 1 & 1\\ 
    1 & 0 & 0 && 1 & 0 & 0 && 0 & 1 & 0 && 1 & 1 & 1\\ 
    1 & 0 & 0 && 0 & 0 & 1 && 0 & 1 & 0 && 1 & 1 & 1\\ 
    0 & 1 & 0 && 1 & 0 & 0 && 0 & 1 & 0 && 1 & 1 & 1\\ 
    0 & 1 & 0 && 0 & 0 & 1 && 0 & 1 & 0 && 1 & 1 & 1\rlap{.}
  \end{tabular}
  \]
  Projecting this $3$-cube down to $\RR^3\times\RR^3$ by $\pi_\Box$ yields, in
  order, $\frac13$~times the points
  \[
  \begin{tabular}[c]{*{7}{c}}
    2 & 0 & 1 && 1 & 1 & 1 \\
    1 & 0 & 2 && 1 & 1 & 1 \\
    1 & 1 & 1 && 1 & 1 & 1 \\
    0 & 1 & 2 && 1 & 1 & 1 \\
    2 & 1 & 0 && 1 & 1 & 1 \\
    1 & 1 & 1 && 1 & 1 & 1 \\
    1 & 2 & 0 && 1 & 1 & 1 \\
    0 & 2 & 1 && 1 & 1 & 1\rlap{,} 
  \end{tabular}
  \]
  and applying $\pi_\Delta$ now has the effect of eliminating the last three
  coordinates. We have found the cell $\pi(\sigma\cap L)$ of
  $\pi\Hom(G_\Gamma,H_\Gamma)$, namely $\frac13$~times the convex hull of the
  points
  \[
     (2,0,1),\;(1,0,2),\;(1,1,1),\;(0,1,2),\;
     (2,1,0),\;(1,1,1),\;(1,2,0),\;(0,2,1),     
  \]
  which is indeed the hexagon $Q_\Gamma(1,1,1)$. Notice how two antipodal
  vertices of~$\sigma\cap L$ are projected down to the same point $(1,1,1)$ in
  the interior of~$Q_\Gamma(1,1,1)$, and thus play no role in the convex hull
  of~$\pi(\sigma\cap L)$. 
\end{example}

\section{The case $G=K_g$}

In this paper, we will focus especially on the case where $G$~is the
complete graph on $g$~vertices, with $V(G)=[g]$.  The case $H=K_h$ has
been widely studied in connection with coloring problems on graphs;
see~\cite{Kozlov05} for a survey.

\subsection{Projections and orbits of the symmetry group}
Both complexes $\Hom_+(K_g,H)$ and $\Hom(K_g,H)$ admit an $S_g$-action
by permuting the vertices of~$K_g$, where $S_g$~is the symmetric group
on $g$~letters. By Proposition~\ref{prop:quotient}, and in contrast to
the situation of Observation~\ref{obs:not_equivalence}, it is now the
case that each cell of~$\pi\Hom(K_g,H)$ represents an
$S_g$-equivalence class of faces of the polytopal
complex~$\Hom(K_g,H)$, so that we can think of $\pi\Hom(K_g,H)$ as a
``quotient'' $\Hom(K_g,H)/S_g$.  Any cell~$\pi(\sigma)$
of~$\pi\Hom(K_g,H)$ is a product of simplices, because
$\sigma_i\cap\sigma_j=\emptyset$ for all $i\ne j$ by the definition of
$\Hom(K_g,H)$ and because $H$~is loopless.  However, in general
$\pi\Hom(K_g,H)$ is not a polytopal complex; we will give a
characterization for when this happens in
Theorem~\ref{thm:hypersimplex} below.  Before this, we prove that
$\pi\Hom_+(K_g,H)$ and $\pi\Hom_+^t(K_g,H)$ have analogous properties
to~$\pi\Hom(K_g,H)$:

\begin{theorem}\label{thm:quotient_kg}\
  \begin{compactenum}[\upshape(i)] \samepage
  \item $\pi_\Box\Hom_+(K_g,H)$ is a simplicial immersion of
    $\Hom_+(K_g,H)$ into $\Delta_{V(H)}\times\Delta_{[g]}$. This means
    that locally on each simplex of~$\Hom_+(K_g,H)$, the projection
    $\pi_\Box$ is a bijection onto a simplex whose vertices are among
    the vertices of~$\Delta_{V(H)}\times\Delta_{[g]}$, but the images
    of different faces may intersect.
    Put differently, $\Hom_+(K_g,H)$ is a ``horizontal'' complex, i.e., it has
    no faces in the kernel of~$\pi_\Box$.
    In summary, the following diagram commutes:
  \[
  \begin{matrix}
    \displaystyle\join_{i\in[g]}\Delta_{V(H)} 
    &\supset&
    \Hom_+(K_g,H) 
    &\xmapsto{\ \iota_L\ }&
    \Hom(K_g,H)\times\bg 
    &\subset&
    \Delta_{[gh]}\times\bg \\[1.5ex]
    &&\downmapsto\pi_\Box && \downmapsto\pi_\Box \\[1.5ex]  
    \Delta_{V(H)}\times\Delta_{[g]}
    &\supset&
    \pi_\Box\Hom_+(K_g,H)
    &\xmapsto{\ \iota_{\pi_\Box L}\ }&
    \pi\Hom(K_g,H)\times\bg
    &\subset&
    \Delta_{V(H)}\times\bg \\[1.5ex]
    &&\downmapsto\pi_\Delta && \downmapsto\pi_\Delta \\[1.5ex]  
    \Delta_{V(H)} 
    &\supset &
    \pi\Hom_+(K_g,H)
    &&
    \pi\Hom(K_g,H)
    &\subset&
    \Delta_{V(H)}
  \end{matrix}
  \]

\item Each cell of $\pi\Hom_+(K_g,H)$ represents an
  $S_g$-equivalence class of faces of the simplicial complex~$\Hom_+(K_g,H)$.
  
  \item The same statements hold with $\Hom_+(K_g,H)$ replaced by
    $\Hom_+^t(K_g,H)$.
  \end{compactenum}
  
\end{theorem}

\begin{proof}
  (i) We first check that $\pi_\Box(\sigma)$~is a simplex. For this,
  note that $\sigma_i\cap\sigma_j=\emptyset$ for any $i\ne j\in[g]$ by
  definition of~$\Hom_+(K_g,H)$, because $K_g$ is complete and $H$~is
  loopless.  Therefore, all simplices $\sigma_i\times e_i$ lie in skew
  subspaces of~$\RR^h\times\RR^g$, and their convex hull is 
  a simplex.  It remains to check that no face~$\sigma$
  of~$\Hom_+(K_g,H)$ lies in the kernel of~$\pi_\Box$.  For this,
  suppose that there is some edge~$\sigma$ in~$\Hom_+(K_g,H)$ that
  gets mapped to a point~$w$ by~$\pi_\Box$. Then
  $\sigma=\star_{i\in[g]}\sigma_i$ with $\sigma_j=\sigma_k=\{w\}$ for
  some $j\ne k$ and $\sigma_i=\emptyset$ for $i\ne j,k$, but this
  again contradicts the fact that $\sigma_i\cap\sigma_j=\emptyset$ for
  any $i\ne j$.
  
  Therefore, the restriction of $\pi_\Box$ to any face of~$\Hom_+(K_g,H)$ is a
  bijection onto some simplex that is the convex hull of vertices
  of~$\Delta_{V(H)}\times\Delta_{[g]}$, but the images of different faces
  of~$\Hom_+(K_g,H)$ may in general intersect.
  
  Part (ii) is now an easy consequence of noting that $\pi\Hom_+(K_g,H)$ is
  indeed a simplicial subcomplex of~$\Delta_{V(H)}$, and that we obtain all
  $g!$~faces in the $S_g$-orbit of a given face~$\sigma\cong\pi_\Box(\sigma)$
  by lifting each simplex $\sigma_i\subset\Delta_{V(H)}$ to $\sigma_i\times
    e_{\pi(i)}$, for all $\pi\in S_g$.

  Finally, (iii) follows from Theorem~\ref{thm:intersect} (i).
\end{proof}

In view of this theorem, $\pi\Hom_+^t(K_g,H)$ is the simplicial subcomplex of
$\Delta_{V(H)}$ that is induced by the family of simplices 
\[
   \{\Delta_V: V \text{ is the
     vertex set of a complete $g$-partite subgraph of~$H$}\}.
\]
The analogous construction $\pi\Hom_+(K_g,H)$ is not very
interesting: since $\Hom_+(K_g,H)$ contains faces of the form
$(\Delta_{V(H)},\emptyset,\dots,\emptyset)$, we just get the whole
simplex~$\Delta_{V(H)}$.

\begin{remark} At the level of faces, the left column of the diagram
  of Theorem~\ref{thm:quotient_kg} reads
  \[
  \begin{matrix}
    \sigma & = & \ \conv\ \bigcup_{i=1}^g\mu_i(\sigma_i)\times e_i \ &
    \subset & \Hom_+^{(t)}(K_g,H)
    & \subset & \displaystyle\join_{i\in[g]}\Delta_{V(H)}\\[2ex]
    &&\downmapsto\pi_\Box && \downmapsto\pi_\Box && \downmapsto\pi_\Box\\[2ex]
    \pi_\Box(\sigma) & = & 
    \conv\ \bigcup_{i=1}^g \sigma_i\times e_i
    & \subset & \pi_\Box\Hom_+^{(t)}(K_g,H)
    & \subset & \Delta_{V(H)}\times\Delta_{[g]}\\[2ex]
    &&\downmapsto\pi_\Delta &&\downmapsto\pi_\Delta &&\downmapsto\pi_\Delta\\[2ex]
    \pi(\sigma) & = & 
    \conv\ \bigcup_{i=1}^g \sigma_i
    & \subset & \pi\Hom_+^{(t)}(K_g,H) & \subset & \Delta_{V(H)}
  \end{matrix}
  \]
\end{remark}

\medskip

\subsection{When is \boldmath $\pi\Hom(K_g,H)$ a polytopal complex?} 

To answer this question, we introduce some special subcomplexes of
$\Hom$-complexes:

\begin{definition}
  Let $G,H$ be loopless graphs. The \emph{induced $\Hom$-complexes}
  $\IHom(G,H)$ and $\IHom_+(G,H)$ are the subcomplexes of~$\Hom(G,H)$,
  respectively $\Hom_+(G,H)$, obtained by considering only induced
  complete bipartite subgraphs.
\end{definition}

Recall that the \emph{clique number} $\omega(H)$ of $H$ is the
number of vertices in a largest clique~in~$H$.

\begin{proposition}\label{prop:ihom}
  Let $H$ be a loopless graph with $\omega(H)\ge g\ge1$. Then
  $\IHom(K_g,H)=\Hom(K_g,H)$ and $\IHom_+(K_g,H)=\Hom_+(K_g,H)$ if and
  only if $\omega(H)=g$.
\end{proposition}

\begin{proof}
  A complete $g$-partite subgraph
  $\sigma=\sigma_1\cup\dots\cup\sigma_g$ of~$H$ is not induced if and
  only if there exists an edge of~$H$ connecting two vertices $x,y$ in
  the same part of $\sigma$, which we may assume to be $\sigma_1$; but
  then we can find a clique of size~$g+1$ by taking $x$,~$y$, and one
  vertex each from $\sigma_2$,~\dots,~$\sigma_g$. Conversely, any
  clique $\{v_1,\dots,v_{g+1}\}$ of size~$g+1$ yields a non-induced
  complete $g$-partite graph
  $\sigma=\{v_1,v_2\}\cup\{v_3\}\cup\dots\cup\{v_{g+1}\}$.
  (Finally, note that $\Hom(K_g,H)=\emptyset$ if $\omega(H)<g$.) 
\end{proof}

Another ingredient to our answer is the \emph{hypersimplex}
$\Delta(h,g)\subset\RR^h$, the convex hull of the set of all
$0/1$-vectors of length $h$ with exactly $g$~ones. It is an
$(h-1)$-dimensional polytope that may be thought of as the slice of
the $h$-dimensional $0/1$-cube with the plane
$\{\xx\in\RR^h:\sum_{i=1}^hx_i=g\}$.  All faces of a hypersimplex are
again hypersimplices, and its edges are spanned by pairs of vertices
whose coordinates differ in exactly two locations.

\begin{proposition}\label{prop:hypersimplex}
  For any (loopless) graph $H$, the $1$-skeleton of $\pi\Hom(K_g,H)$
  is a subcomplex of the $1$-skeleton of the
  hypersimplex~$\Delta(h,g)$.
\end{proposition}

\begin{proof}  
  By abuse of notation, we replace $\pi\Hom(K_g,H)$ by a $g$ times inflated
  copy but keep the same name. Let $\sigma$ be a face of $\Hom(K_g,H)$, so
  that $\sigma_i\cap\sigma_j=\emptyset$ for all $1\le i\ne j\le g$. As in the
  proof of Proposition~\ref{prop:quotient}, we can write any point
  $z\in\pi(\sigma)$ as
  \[
     z 
       \ = \  \sum_{i\in V(G)} \sum_{v\in\vertices\sigma_i}\lambda_{iv}v\,.
  \]
  The vertices of $\pi(\sigma)$ are $0/1$-vectors of length $h$ (because the
  vertices of the $\sigma_i$'s are vertices of $\Delta_H\subset\RR^h$) with
  exactly $g$ ones (there are $g$ mutually disjoint simplices~$\sigma_i$, and
  to get a vertex of~$\sigma$, exactly one $\lambda_{iv}$ must be~$1$, for
  each $i$); in other words, they are vertices of~$\Delta(h,g)$.  The
  statement now follows because the sets $\sigma_i$ of any edge~$\sigma$
  of~$\Hom(K_g,H)$ all have size~$1$, except for exactly one set of size~$2$.
  Therefore, the coordinates of the vertices of $\pi(\sigma)$ differ in
  exactly two places.
\end{proof}

We are now in a position to answer the question posed at the beginning of this
section:

\begin{theorem}\label{thm:hypersimplex}
  $\pi\Hom(K_g,H)$ is a non-empty polytopal complex, and
  $\pi\Hom_+(K_g,H)$ a non-empty simplicial complex, if and only if
  $\omega(H)=g$.
\end{theorem}

\begin{proof}
  Let $\supp:\RR^h\to[h]$ denote the map that assigns to any vector
  $v\in\RR^h$
  the set of all indices $i\in[h]$ such that $v_i\ne0$. 
  Moreover, for any
  face $\sigma=(\sigma_1,\dots,\sigma_g)$ of~$\Hom(K_g,H)$, let $G_\sigma$
  denote the complete $g$-partite subgraph of~$H$ on the vertex set
  $V_\sigma=\sigma_1\cup\dots\cup\sigma_g\subset[h]$ that is associated
  to~$\sigma$ by the definition of faces of~$\Hom$-complexes. 
  
  We want to rule out that two faces $\pi(\sigma)$ and $\pi(\tau)$
  of~$\pi\Hom(K_g,H)$ intersect badly, where $\sigma,\tau$~are faces
  of~$\Hom(K_g,H)$. If they do, there exists a circuit (i.e., a minimal affine
  dependency) $C=C_+\cup C_-$ among the vertices of~$\Delta(h,g)$, such that
  $C_+\subset\vertices\pi(\sigma)$, $C_-\subset\vertices\pi(\tau)$ and
  $C_+\cap C_-=\emptyset$.  In this case, we can find sets of
  positive real numbers $\{\lambda_v:v\in C_+\}$ and $\{\mu_w:w\in C_-\}$ with
  $\sum_{v\in C_+}\lambda_v=1$ and $\sum_{w\in C_-}\lambda_w=1$, such that
  $\sum_{v\in C_+}\lambda_v v = \sum_{w\in C_-}\lambda_w w$.  In particular,
  \begin{equation}\label{eq:supp}
    \supp\sum_{v\in C_+}v \ = \ \supp\sum_{w\in C_-}w \ =: \
    V\subset[h]\,.
  \end{equation}
  We claim that $V$~is the vertex set of an ---~at this point not necessarily
  unique~--- complete $g$-partite subgraph of~$H$ of the form~$G_\rho$ for
  some face~$\rho$ of~$\Hom(K_g,H)$. Indeed, each vertex in~$C$ yields a
  $g$-clique in~$H$ contained in the vertex set~$V$. Because $C_+\subset
  \vertices\pi(\sigma)$, the vertices in~$C_+$ yield a complete (but not
  necessarily induced) $g$-partite graph~$G_{\rho_+}$ on the vertex set~$V$
  that corresponds to a face $\rho_+$ of~$\sigma$.  The same happens
  for~$C_-\subset\vertices\pi(\tau)$, and we obtain another complete
  $g$-partite graph~$G_{\rho_-}$ on the \emph{same} vertex set~$V$.
  
  If all complete $g$-partite subgraphs of~$H$ are induced, then
  $G_{\rho_+}=G_{\rho_-}$; otherwise, some edge of~$G_{\rho,-}$, say, would
  join two vertices belonging to the same part of~$G_{\rho,+}$, and the latter
  graph would not be induced. Therefore, the
  graph~$G_\rho:=G_{\rho_+}=G_{\rho_-}$ is unique, and with it $V:=V_\rho$ and
  the face $\rho:=\rho_+=\rho_-$ of~$\Hom(K_g,H)$. But then $\rho$~is a common
  face of both $\sigma$~and~$\tau$, and~\eqref{eq:supp} says that the
  circuit~$C$ is supported on the common face~$\pi(\rho)$ of
  $\pi(\sigma)$~and~$\pi(\tau)$ (which is a product of simplices by
  Theorem~\ref{thm:permutohedron}). In consequence,
  $\pi(\sigma)$~and~$\pi(\tau)$ do not intersect badly after all.
  
  Conversely, if $G_{\rho_+}\ne G_{\rho_-}$ are different complete $g$-partite
  graphs on the common vertex set~$V$ given by~\eqref{eq:supp}, then the face
  $\rho_-$ of~$\tau$ is not a face of~$\sigma$, and the projections
  $\pi(\sigma)$~and~$\pi(\tau)$ do have a bad intersection.
  The theorem now follows from Proposition~\ref{prop:ihom}.
  \end{proof}

\begin{example}
  $\pi\Hom(K_2,K_4)$ is not a polytopal complex, as predicted by
  Theorem~\ref{thm:hypersimplex}: not all complete bipartite graphs
  in~$K_4$ are induced. The complex $\Hom(K_2,K_4)$ is isomorphic to the
  boundary complex of the $3$-dimensional cuboctahedron, and by
  Proposition~\ref{prop:hypersimplex}, the $1$-skeleton of $\pi\Hom(K_2,K_4)$
  is contained in the $1$-skeleton of the octahedron~$\Delta(4,2)$ (they
  actually coincide in this case). However, the six square faces
  of~$\Hom(K_2,K_4)$ get mapped to three square faces in~$\pi\Hom(K_2,K_4)$,
  namely the three ``internal squares'' of the octahedron~$\Delta(4,2)$, and
  the relative interiors of any two of these intersect.
\end{example}

\begin{remark}
  The vertices of $\pi\Hom(K_g,H)$ are those vertices of the hypersimplex
  $\Delta(h,g)$ that correspond to cliques in~$H$, and by
  Proposition~\ref{prop:hypersimplex} the $1$-skeleton of each
  cell~$\pi(\sigma)$ is entirely contained in the $1$-skeleton
  of~$\Delta(h,g)$. Thus, each cell~$\pi(\sigma)$ is a \emph{matroid
    polytope}~\cite{Kapranov93}, corresponding to the associated \emph{clique
    matroid} $M_\sigma$. The elements of~$M_\sigma$ are the vertices of the
  complete $g$-partite subgraph~$\sigma$ of~$H$, and its bases the $g$-cliques
  of~$\sigma$. These matters, as well as the connections to tropical geometry,
  are however beyond the scope of this article and will be pursued in a future
  publication. 
\end{remark}

\begin{remark}
  We have seen that $\pi\Hom(K_g,H)\times\bg =
  \pi_\Box\Hom_+(K_g,H)\cap(\RR^h\times\bg)$ is the ``slice parallel
  to~$\Delta_{[g]}$'' given by the polyhedral Cayley trick. Therefore
  it seems natural to ask about
  $\Sigma_h:=h\cdot\pi_\Box\big(\Hom_+(K_g,H)\big) \cap
  (\boldsymbol{\frac{1}{h}}\times\RR^g)$, the ``slice parallel
  to~$\Delta_{V(H)}$''. If $\sigma=\sigma_1\join\cdots\join\sigma_g$ is a cell
  of $\Hom_+(K_g,H)$, then as before any point
  $x\in\pi_\Box(\sigma)=\bigcup_{i=1}^g\sigma_i\times e_i$ can be
  written as a convex combination 
  \[
      x \ = \ \sum_{i=1}^g\lambda_i 
      \sum_{v\in\vertices \sigma_i}\lambda_{i,v} v\times e_i,
  \]
  so that the intersection of $\pi_\Box(\sigma)$ with
  $\boldsymbol{\frac{1}{h}}\times\RR^g$ is non-empty if and only if
  $\bigcup_{i=1}^g\sigma_i = [h]$, where we view the $\sigma_i$ as subsets
  of~$[h]$. In other words, every complete $g'$-partite subgraph of~$H$
  supported on \emph{all} $h$~vertices of~$H$ contributes a cell to
  $\Sigma_h$, for all $1\le g'\le g$. For instance, the cells
  of~$\Hom_+(K_g,H)$ of the form
  $\Delta_{V(H)}\join\emptyset\join\dots\join\emptyset$ each contribute a
  vertex of the form  $b\times e_i$, where $b$~is the barycenter
  of~$\Delta_{V(H)}$. 
\end{remark}

\section{Dissection complexes}

For $k\ge3$ and $m\ge1$, consider the set of dissections of a convex
$N$-gon into $m$~convex $k$-gons. Note that for such a dissection to
be possible, it is necessary and sufficient that $N=m(k-2)+2$. We
agree to label the vertices of the $N$-gon by $0$, $1$, \dots, $N-1$.
Let $\delta(k,m)$ be the set of \emph{$k$-allowable} diagonals of the
$N$-gon, i.e., those diagonals that can appear in a dissection into
$k$-gons. It is easy to check that these are precisely the diagonals
that connect a vertex~$x$ with one of the form $x+k-1+j(k-2)\bmod N$,
for $0\le j\le m-2$, and that $|\delta(k,m)|=(m-1)N/2$.
Let $\cross(\delta(k,m))$ and $I(k,m)=\Ind(\delta(k,m))$ be the
\emph{crossing graph} and \emph{independence graph}
of~$\delta(k,m)$. These are complementary graphs on the vertex
set~$\delta(k,m)$, such that two vertices are joined by en edge in
$\cross(\delta(k,m))$ if the corresponding $k$-allowable diagonals
intersect in their relative interior, while the same two vertices are
joined in~$I(k,m)$ if this is not the case. For a graph $G$, the
simplicial complexes $\IInd(G)$ and $\Cl(G)$ are the
\emph{independence complex} and \emph{clique complex} of~$G$, whose
simplices are the independent sets, respectively the cliques,
of~$G$. Thus, $I(k,m)=\sk^1\IInd(\delta(k,m))$ and
$\cross(\delta(k,m))=\sk^1\Cl(\delta(k,m))$.

\begin{proposition}\label{prop:ICDelta}
  $\Hom\big(K_{m-1}, I(k,m)\big)$ is the polytopal complex
  on the vertex set $\delta(k,m)$ whose cells are the products
  $\Delta^{C_1}\times\dots\times\Delta^{C_{m-1}}$ of $m-1$ simplices such that
  each $C_i$ is a non-empty clique in~$\cross(\delta(k,m))$, and
  $C_i$~and~$C_j$ are independent in $\cross(\delta(k,m))$ for $i\ne j$.
\end{proposition}

\begin{proof}
  $\Hom\big(K_{m-1}, I(k,m)\big)$ 
  arises by labelling the vertices of $K_{m-1}$ with non-empty lists
  of elements in $\delta(k,m)$, such that any two diagonals of two
  lists $\lambda$, $\lambda'$ on different vertices are joined by an
  edge in $I(k,m)$. But this means precisely that the diagonals in
  $\lambda$ do not cross the diagonals in $\lambda'$. Moreover, no
  list can contain two or more independent diagonals, because on the
  one hand all $m-1$ lists must be non-empty, and on the other hand
  the maximal size of an independent set in $\delta(k,m)$ is~$m-1$ by
  definition.
\end{proof}

\begin{remark}
  That $\widetilde\DD(k,m)=\Hom\big(K_{m-1}, I(k,m)\big)$ is a polytopal
  complex follows from its definition as $\Hom$-complex, but can also be
  proven directly from Proposition~\ref{prop:ICDelta}:
  First note that any face of a cell of~$\widetilde\DD(k,m)$ is again a
  product of $m-1$~simplices, and therefore again a cell
  of~$\widetilde\DD(k,m)$.  On the other hand, let $C=\prod_{i=1}^{m-1}
  \Delta^{C_i}$ and $D=\prod_{i=1}^{m-1}\Delta^{D_i}$ be two cells
  of~$\widetilde\DD(k,m)$.  We must show that $C\cap D$ is either a common
  face of both $C$ and $D$, or else does not index a face in the complex. For
  this, let $D_{\pi(1)}$, \dots, $D_{\pi(m-1)}$ for $\pi\in S_{m-1}$ be a
  permutation of the $D_i$ such that $C_i\cap D_{\pi(i)}\ne\emptyset$ for all
  $1\le i\le m-1$. If such a permutation~$\pi$ does not exist, then the
  intersection $\vertices C \cap \vertices D$ is not a union of $m-1$
  independent cliques, and so $C\cap D$~is~not a face of~$\widetilde\DD(k,m)$.
  On the other hand, if such a permutation exists, then it is unique: Suppose
  that $C_i\cap D_{\pi(i)}=V_i\ne\emptyset$ and $C_i\cap
  D_{\sigma(i)}=W_i\ne\emptyset$ for some $1\le i \le m-1$ and $\pi,\sigma\in
  S_{m-1}$.  Then $V_i,W_i\subset C_i$, so that $V_i$ and $W_i$ are not
  independent; but then $V_i\subset D_{\pi(i)}$, $W_i\subset D_{\sigma(i)}$
  forces $\pi(i)=\sigma(i)$.  The cell of $\widetilde\DD(k,m)$ corresponding
  to $C\cap D$ is then $\prod_{i=1}^{m-1}\Delta^{C_i\cap D_{\pi(i)}}$.
\end{remark}

The number of dissections of a convex polygon into $m$~convex $k$-gons was
already determined in~1791 by Fuss~\cite{Fuss91}; see also the simplified
proof in~\cite{Przytycki-Sikora00}. Two related complexes were considered
somewhat later: In~2005, Tzanaki~\cite{Tzanaki05} proved that the simplicial
complex $\TT(k,m)=\IInd(\delta(k,m))$ is homotopy equivalent to a wedge of
$\frac{1}{m}\binom{m(k-2)}{m-1}$ spheres of dimension $m-2$. Its dual graph is
the \emph{flip graph} $D(m,k)$, whose vertices are the dissections of the
polygon (the facets of~$\TT(k,m)$), and in which two dissections are adjacent
if they differ in the placement of exactly one interior
diagonal~\cite{Huemer-Hurtado-Pfeifle05}.

\begin{definition} We will use the following abbreviations:
  \begin{eqnarray*}
    \DD(k,m) &=& \pi\Hom\big(K_{m-1}, I(k,m)\big), \\
    \DD_+(k,m) &=& \pi\Hom_+\big(K_{m-1}, I(k,m)\big),\\
    \DD_+^t(k,m) &=& \pi\Hom_+^t\big(K_{m-1}, I(k,m)\big).
  \end{eqnarray*}
\end{definition}

\begin{proposition}\label{prop:hom}
  With these notations,
  \begin{compactenum}[\upshape(a)]
  \item $\DD(k,m)$ is a polytopal complex, and $\DD_+(k,m)$ and
    $\DD_+^t(k,m)$ are simplicial complexes.  Moreover, $\DD(k,m)$ arises
    as a linear section of $\DD_+(k,m)$ and $\DD_+^t(k,m)$. 
    
  \item\label{prop:hom:t} $\TT(k,m)$ is the simplicial complex induced on the
    set of transversal $(m-2)$-dimensional faces of~$\DD_+(k,m)$.
    
  \item $D(k,m)=\sk^1\DD(k,m)$.
  \end{compactenum}
\end{proposition}

\begin{proof}
  The first assertion of (a) follows from Theorem~\ref{thm:hypersimplex} and
  the definition of~$I(k,m)$, the second one from Theorem~\ref{thm:intersect},
  and the third is true by definition.  For (b), note that the faces of a
  transversal $(m-2)$-dimensional simplex of~$\Hom_+(K_{m-1},I(k,m))$ are
  obtained by labelling each vertex of $K_{m-1}$ with a list of diagonals of
  size~$0$~or~$1$. By the definition of~$I(k,m)$ these diagonals are mutually
  non-crossing, so after dividing out by the $S_{m-1}$-symmetry we obtain
  exactly the simplices of~$\TT(k,m)$.  The reasoning for (c) is similar.
\end{proof}

\begin{figure}[htbp]
  \centering
  \includegraphics[width=.9\linewidth]{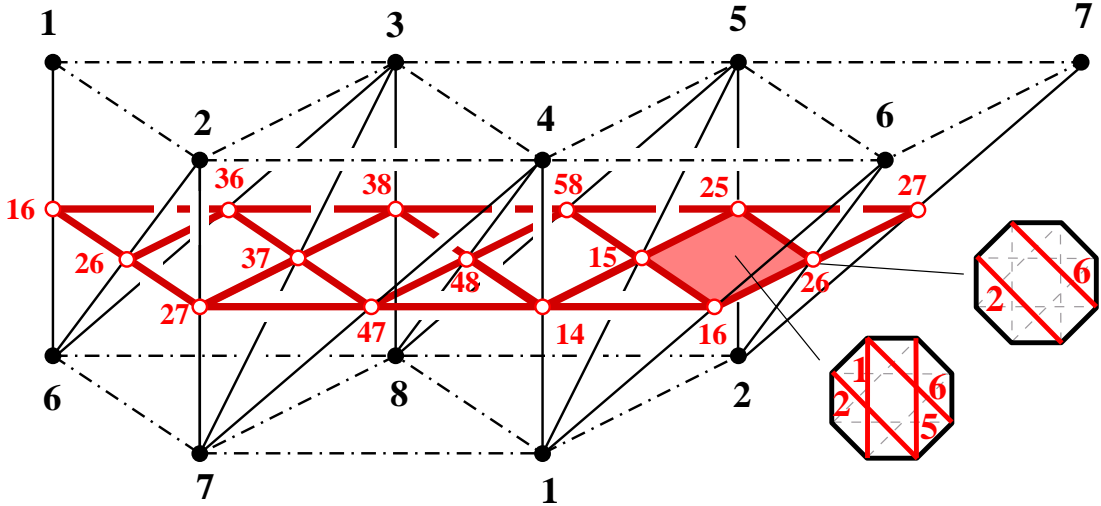}
  \caption{The projected simplicial complex $\DD^t_+(4,3)$ with the
    polytopal complex $\DD(4,3)$ arising as a section.  Notice that $\DD(4,3)$
    is a M\"obius band, while its double cover $\Hom(K_2,I(4,3))$ is
    homeomorphic to $S^1\times[0,1]$. The solid edges of $\DD^t_+(4,3)$ make
    up the complex $\TT(4,3)$. The $4$-allowable diagonals of the
    corresponding $8$-gon are numbered cyclically from~$1$~to~$8$.}
  \label{fig:d43}
\end{figure}

\begin{remark}
  $\DD_+(k,m)$ is not the only simplicial complex that contains $\TT(k,m)$ as
  a subcomplex. For example, one can also consider the simplical complex
  $\IC_\Delta(k,m)$ on the vertex set~$\delta(k,m)$ whose simplices are unions
  of independent cliques in~$\cross(\delta(k,m))$. This means that all
  diagonals in any of the cliques intersect in their relative interior, but
  any two diagonals from different cliques do not.  $\IC_\Delta(k,m)$ is a
  strict subcomplex of~$\DD_+(k,m)$, because any tuple of mutually independent
  cliques is a tuple of mutually independent sets of diagonals, but not vice
  versa.  Moreover, one can check that $\IC^t_\Delta(k,m)=\DD_+^t(k,m)$, where
  $\IC^t_\Delta(k,m)$ is the induced simplicial complex on the set of unions
  of $m-1$~non-empty mutually independent cliques. Therefore, $\DD(k,m)$ is
  also a linear section of $\IC_\Delta(k,m)$.
\end{remark}

\subsection{On the dimension and homotopy type of the dissection complexes}

\begin{proposition}
  $\dim\DD(k,m)=\big\lfloor\frac{m}{2}\big\rfloor(k-2)$, and
  $\dim\DD_+(k,m)=\dim\DD(k,m)+m-2$.
\end{proposition}

\begin{proof}
  The second statement holds generally, because the dimension
  of a maximal cell $\sigma_1\times\dots\times\sigma_{m-1}$ is
  $|\sigma_1|+\dots+|\sigma_{m-1}|-m+1$, while
  $\dim\sigma_1\star\dots\star\sigma_{m-1}=
  |\sigma_1|+\dots+|\sigma_{m-1}|-1$.
  
  To prove the first statement, we must find $m-1$ sets of mutually
  intersecting $k$-allowable diagonals of the $N$-gon in such a way that no
  two diagonals from distinct sets cross, and the total number of diagonals is
  maximized. For this, we fix diagonals $d_1$, \dots, $d_{m-1}$ of the
  $N$-gon. To each $d_i$ (with endpoints $x_i$~and~$y_i$) and some choice of
  $\delta_i^+,\delta_i^-\in\NN$, we adjoin all $k$-allowable diagonals with one
  endpoint between $x_i-\delta_i^-$ and $x_i+\delta_i^+$, and the other
  endpoint between $y_i-\delta_i^+$ and $y_i+\delta_i^-$
  (cf.~Figure~\ref{fig:dimension}).
  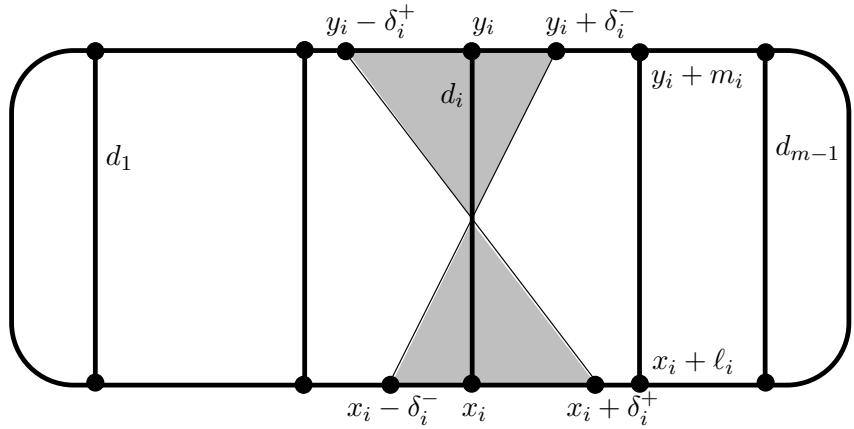
\begin{figure}[htbp]
    \centering
    \input{dimension.pstex_t}    
    \caption{Finding a cell of $\DD(k,m)$ of maximal dimension.}
    \label{fig:dimension}
  \end{figure}
  Note that $0\le\delta_i^+ + \delta_{i+1}^- \le \ell_i$ and $0\le
  \delta_i^-+\delta_{i+1}^+ \le m_i$ for $1\le i \le m-2$, where
  $\ell_i=x_{i+1}-x_i$ and $m_i=y_{i+1}-y_i$, so that $\ell_i+m_i=k-2$. All
  calculations are to be considered modulo~$N$. The dimension of the
  cell~$\sigma=\sigma_1\times\dots\times\sigma_{m-1}$ obtained in this way is
  then
  \begin{eqnarray*}
     -m+1+\sum_{i=1}^{m-1} |\sigma_i| & = & 
     \sum_{i=1}^{m-1} \delta_i^-+\delta_i^+ \\ &\le & 
     \min\left\{\delta_1^- +\delta_{m-1}^+ + \sum_{i=1}^{m-2}\ell_i,\quad
       \delta_1^+ +\delta_{m-1}^- + \sum_{i=1}^{m-2}(k-2-\ell_i)\right\}.
  \end{eqnarray*}
  For even~$m$, we can maximize this value by choosing
  $\ell_{2j-1}=\big\lfloor\frac{k-2}{2}\big\rfloor$ and
  $\ell_{2j}=\big\lceil\frac{k-2}{2}\big\rceil$ for $1\le j\le
  \frac{m-2}{2}$; moreover, it is readily verified that in this case
  $\delta_1^-$~and~$\delta_{m-1}^+$ can be chosen such that
  $\delta_1^-+\delta_{m-1}^+ = k-2$. Therefore,
  \[
     \dim\sigma \ = \ k-2 + (m-2)\frac{k-2}{2} \ = \ \frac{m}{2}(k-2),
  \]
  as claimed. The proof for odd~$m$ is similar.
\end{proof}

The simplicial complex $\DD_+(k,m)$ and the polytopal one $\DD(k,m)$ are
homotopy equivalent, because the former arises by taking as simplices the
union of vertices of polytopal cells of the latter; cf.~\cite{Dong02}.  Since
$\DD_+(k,m)$~is obtained from~$\TT(k,m)$ by adjoining extra cells, and
$\TT(k,m)$ is homotopy equivalent to a wedge of ``many'' $(m-2)$-dimensional
spheres, it seems reasonable to hope for the extra cells to simplify the
topology.

Explicit computation in small cases reveals that this is indeed the case, but
perhaps not to the greatest extent possible. Although the homology of
$\DD_+(k,m)$, and therefore of~$\DD(k,m)$, is substantially simpler than that
of $\TT(k,m)$, Table~\ref{tab:homology} below does not support the conjecture
that these complexes are homotopy equivalent to a wedge of a ``simple'' number
of spheres. 

\begin{table}[htbp]
  \centering\small
  \begin{tabular}[c]{c|*{7}{c}}
    $k\setminus m$ & 3 & 4 & 5 & 6 & 7 & 8 & 9\\\hline
    \raisebox{-1ex}{3} & \rule{0pt}{2.7ex} 2 (1) & 3--4 (1--2) & 4--5 (1--2) &
    5--7 (1--3) & 7--8 (2--3) &  8--10 (2--4) & 9--11 (2--4)\\ 
    & $r_1=1$ & $r_2=1$ & $r_3=1$ & $r_4=1$ & $r_5=1$ & $r_6=1$ \\[2ex]
    \smash{\raisebox{-2.5ex}{4}} & 3 (2) & 4--6 (2--4) & 5--7 (2--4) & 6--10
    (2--6) & 7--11 (2--6) & 8--14\\ 
    & $r_1=1$ & $r_3=1$ & \parbox[t]{1.5cm}{\centering $r_3=4$
    $r_4=4$} & $r_5=1$ &  \parbox[t]{1.5cm}{\centering $r_5=17$
    $r_6=20$}\\[4ex]  
    \raisebox{-1ex}{5} & 4 (3) & 5--8 (3--6) & 6--9 (3--6) & 7--13 (3--9) &
    8--14\\ 
    & $r_1=1$ & $r_3=1$ & $r_3=1$ & $r_5=17$\\[2ex]
    \raisebox{-1ex}{6} & 5 (4) & 6--10 (4--8) & 7--11 (4--8) & 8--16\\
    & $r_1=1$ & $r_3=1$ & $r_3=1$\\[2ex]
    \raisebox{-1ex}{7} & 6 (5) & 7--12 (5--10) & 8--13\\
    & $r_1=1$ & $r_3=1$
  \end{tabular}
\bigskip
  \caption{Dimensions of facets and nonzero integer homology ranks of 
    $\DD_+(k,m)$.  The first line of each
    entry~$(k,m)$ 
    lists the range of dimensions of the facets of $\DD_+(k,m)$, respectively
    $\DD(k,m)$,
    and the next lines the nonzero ranks~$r_i$ of their reduced integer
    homology 
    groups, so that $\tilde H_i(\DD_+(k,m), \ZZ)=\ZZ^{r_i}$. No
    homology groups in the table have torsion.
  } \label{tab:homology}
  
\end{table}

\section{Staircase triangulations and cyclic polytopes}
\label{sec:staircase}

In Proposition~\ref{prop:hom}, we found $\TT(k,m)$ and $D(k,m)$ as
subcomplexes of $\DD_+(k,m)$, respectively of~$\DD(k,m)$.  Next, we
identify other ``nice'' subcomplexes of~$\DD(k,m)$ and its relatives:

\begin{theorem}\label{thm:subcomplexes}
  Let $r,s$ be integers such that $1\le r\le
  \frac{(m-1)(k-2)+2}{k-1}$ and $1\le s \le k-1$.

  \begin{compactenum}[\upshape(a)]
  \item The simplicial complex $\DD_+(k,m)$ contains copies of the staircase
    triangulation~$\Sigma(r,s)$ of the product of simplices
    $\Delta^{r-1}\times\Delta^{s-1}$.
    
  \item The polytopal complex $\DD(k,m)$ contains copies of the polytopal
    complex $\CC(r,s)$, where $d=2s-2$ and $n=r+d$.
  \end{compactenum}

\end{theorem}

In order to define the complexes mentioned in this theorem, we must
first recall some facts about cyclic polytopes. A standard realization
of the \emph{cyclic polytope} $C_d(n)\subset\RR^d$ is given by the
convex hull of any $n$~distinct points $\mu(t_1),\dots,\mu(t_n)$ on
the \emph{moment curve} $\mu:\RR\to\RR^d, t\mapsto(t,t^2,\dots,t^d)$,
where we assume $t_1<\dots<t_n$. Implicit in this definition is the
fact that the combinatorial type of~$C_d(n)$ does not depend on the
concrete values of the $t_i$. A set $I\subset[n]$ indexes a face of
$C_d(n)$ if $I$~satisfies \emph{Gale's evenness criterion}: For any
$j,k\in[n]\setminus I$, the number of elements of~$I$ between
$j$~and~$k$ must be even. Henceforth, we will always identify faces of
$C_d(n)$ with their index sets. A facet of~$C_d(n)$ indexed by
$I\subset[n]$ is a \emph{lower facet} if the cardinality of the
end-set of~$I$ is even, where the \emph{end-set} of~$I$ consists of
the last block of contiguous elements of~$I$. We leave it to the
reader to check (or consult in the literature) the fact that the last
entry of the normal vector of any lower facet of a standard
realization of~$C_d(n)$ is negative.

\subsection{Weak compositions and cyclic
  polytopes~\cite{Huemer-Hurtado-Pfeifle05}} Let $r,s\ge 1$ be integers. A
\emph{(weak) composition\footnote{This seems to be the standard name in the
    literature for ordered partitions.} of $r$ into $s$ parts} is an ordered
$s$-tuple $(a_1,a_2,\dots,a_s)$ of non-negative integers such that
$a_1+a_2+\dots+a_s=r$.  We make the set~$C(r,s)$ of all compositions of~$r$
into $s$~parts into a graph by declaring two of them to be adjacent if they
differ by one in exactly two positions that are connected by a (perhaps empty)
sequence of~$0$'s.  For example, the composition $(1,0,2,4,0,1)$ is adjacent
to $(1,0,2,3,0,2),$ but not to $(0,0,2,4,0,2)$.

\begin{definition}
  Let $n\ge d\ge2$, let $d$ be even and $C_d(n)$ be a $d$-dimensional cyclic
  polytope with $n$~vertices. For any set $I=I(F)=\{i_1,i_1+1, \dots, i_{d/2},
  i_{d/2}+1\}\subset[n]$ that indexes a lower facet~$F=F(I)$ of~$C_d(n)$, the
  sequence~$\chi(F)=\chi(I)$ records the sizes of the ``holes'' in~$I$. More
  precisely, $\chi(F):=\big(i_{j+1}-i_j-2:0\le j \le
  d/2\big)$, where $i_0:=-1$ and $i_{d/2+1}:=n+1$.
\end{definition}

For example, if $n=6$, $d=4$ and $I=I(F)=\{2,3,5,6\}$, then $\chi(I)=(1,1,0)$.

\begin{proposition}\label{prop1}
  Let $r,s\ge 1$ be integers and set $d=2s-2$ and $n=r+d$.
  \begin{enumerate}[\hspace{\parindent}\upshape(a)]
  \item\label{prop1:a} $\chi(F)\in\{0,1,\dots,n-d\}^{\frac{d}{2}+1}$ and \ 
    $\sum_{k=1}^{\frac{d}{2}+1} \chi(F)_k = n-d$ \ for any lower facet~$F$
    of~$C_d(n)$.
  \item\label{prop1:b} The map $\chi$ induces a bijection between the set of
    lower facets of $C_{d}(n)$ and the set of vertices of $C(r,s)$ that takes
    a facet $F$ with $\chi(F)=(a_1,\dots,a_s)$ to the weak composition
    $r=a_1+\dots+a_s$.
  \end{enumerate}
\end{proposition}

\begin{proof}
  Part~\eqref{prop1:a} and the forward direction of part~\eqref{prop1:b}
  follow because any facet~$F$ leaves $n-d=r$ ``holes'' in $\{1,2,\dots,n\}$.
  For the other direction of~\eqref{prop1:b}, Gale's evenness criterion
  uniquely reconstructs~$F$ from any weak composition $r=a_1+\dots+a_s$ by
  inserting a pair of indices between each pair of ``holes'' of sizes $a_i$
  and $a_{i+1}$.  As a check, note that the number of lower facets of the
  cyclic polytope~$C_d(r+d)$ is $\binom{r+d-\lceil d/2\rceil}{\lfloor
    d/2\rfloor}=\binom{r+s-1}{s-1}= |C(r,s)|$.
\end{proof}

\subsection{Cyclic polytopes and staircase triangulations}

Let $r,s\ge 1$ be integers.  A \emph{lattice path} in the grid $[r]\times[s]$
is a connected chain of horizontal and vertical line segments of unit length
that connects $(1,1)$ to $(r,s)$ and is weakly monotone with respect to both
coordinates. Thus, any lattice path has $r+s-1$ vertices. In this paper, we
will always think of a lattice path as its set of vertices.
  
  Denote by $\LL(r,s)$ the set of all \emph{partial lattice paths} in
  $[r]\times[s]$, i.e., all subsets of lattice paths. By identifying
  any partial lattice path with its vertex set, we make $\LL(r,s)$
  into a simplicial complex. This simplicial complex also appears in
  the guise of the \emph{staircase triangulation}~$\Sigma(r,s)$ of the
  product of simplices $\Delta^{r-1}\times\Delta^{s-1}$: it is
  straightforward to check that each partial lattice path in fact
  indexes a simplex in this product polytope, and that all these
  simplices combine to a triangulation. 
  
  For any partial lattice path $\lambda\in\LL(r,s)$, let
  $\lambda_j=\lambda\cap(\{j\}\times[s])$ be the $j$-th ``vertical
  slice'' of~$\lambda$, for $j=1,\dots, r$, and let
  $\aa(\lambda)=(a_1,\dots,a_s)$ be the vector whose $i$-th entry is
  $a_i=|\lambda\cap([r]\times\{i\})|$, the cardinality of the $i$-th
  ``horizontal slice''. Now let $\LL^t(r,s)$ be the set of partial
  lattice paths~$\lambda$ such that $|\lambda_j|\ge1$ for all
  $j=1,\dots,r$, and denote the set of inclusion-minimal members of
  $\LL^t(r,s)$ ---~those with $|\lambda_j|=1$ for all~$j$~---
  by~$M(r,s)$.
  
 Define the $(s-1)$-dimensional polytopal complex
  $\CC(r,s)=\Sigma(r,s)\cap L$ to be the intersection of the staircase
  triangulation~$\Sigma(r,s)$ of~$\Delta^{r-1}\times\Delta^{s-1}$ and the
  $(s-1)$-dimensional plane
  $L=(\frac{1}{r},\dots,\frac{1}{r})\times\RR^{s-1}$. The complex $\CC(r,s)$
  is a mixed polyhedral subdivision of the interior of an $r$~times inflated
  standard $(s-1)$-dimensional simplex~$r\Delta^{s-1}$, such that the vertices
  of~$\CC(r,s)$ are precisely the lattice points of~$r\Delta^{s-1}$.  Finally,
  let $\Sigma^t(r,s)$ denote the set of \emph{transversal} simplices of
  $\Sigma(r,s)$, i.e.\ those that have non-empty intersection with~$L$.
  
 Clearly, the simplices in $\Sigma^t(r,s)$ correspond on
  the one hand bijectively to the partial lattice paths
  in~$\LL^t(r,s)$, and on the other hand ---~via the polytopal Cayley
  trick~--- to the polytopal cells of~$\CC(r,s)$. Under this
  bijection, the set $M(r,s)$ of minimal partial lattice paths
  corresponds to the vertices of~$\CC(r,s)$.



\begin{example}\label{ex:c23}
  Figure~\ref{fig:cayley1} (left) illustrates this correspondence for
  $r=2$, $s=3$. In fact, the bottom part of that figure shows
  $C(2,3)$, the graph of (weak) compositions of $2$ into $3$
  non-negative summands, embedded as the $1$-skeleton of a polyhedral
  decomposition of the twice dilated standard simplex~$2\Delta^2$. In
  the top part, we can see~$C(2,3)$ as a section of the transversal
  part~$\Sigma^t(2,3)$ of the staircase triangulation
  of~$\Delta^2\times\Delta^1$.
\end{example}

\begin{theorem}\label{thm:complexes1}
  Let $r,s\ge1$ be integers, and set $d=2s-2$ and $n=r+d$. Then $\CC(r,s)$~is
  isomorphic to the polytopal subcomplex of the polar-to-cyclic polytope
  $C_d(n)^\Delta$ induced on the vertices dual to lower facets, whose faces
  are of dimension at most~$d/2$. Moreover, the $1$-skeleton of each of these
  complexes is the composition graph~$C(r,s)$.
\end{theorem}

\begin{figure}[htbp]
  \centering
  \input{cyclic.pstex_t}
  \caption{An example for the correspondences of Theorem~\ref{thm:complexes1}.}
  \label{fig:cyclic}
\end{figure}
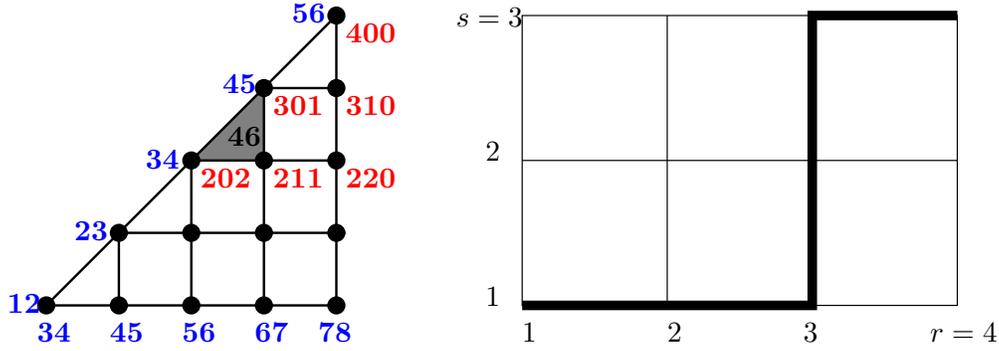

\begin{example}
  Figure~\ref{fig:cyclic} illustrates Theorem~\ref{thm:complexes1} for $r=4$
  and $s=3$ (so that $d=4$ and $n=8$). On the left, we see the
  $s-1=2$-dimensional polytopal complex~$\CC(4,3)$, whose graph is the
  composition graph~$C(4,3)$ of weak compositions of~$4$ into $3$~summands;
  cf.~the labels to the lower right of each vertex.  At the same time, this
  complex is the subcomplex of the polar-to-cyclic polytope~$C_4(8)^\Delta$
  that is induced on the vertices dual to lower facets of~$C_4(8)$: each
  vertex is dual to the lower facet indexed by the union of the labels below
  and to the left of it.
  
  On the right, we see a (full) lattice path in $\LL(4,3)$ that is the
  union of the minimal transversal lattice paths $\{11,21,31,43\}$,
  $\{11,21,32,43\}$ and $\{11,21,33,43\}$ in $M(4,3)$. Their
  $\aa$-vectors are $(3,0,1)$, $(2,1,1)$ and $(2,0,2)$, respectively,
  and so they correspond to the lower facets 4567, 3467~and~3456
  of~$C_4(8)$. The entire lattice path thus corresponds to the face
  of~$\CC(4,3)$ dual to the intersection of these facets, namely the
  $d/2=2$-dimensional triangular face dual to the edge~46.
\end{example}

\begin{proof}[Proof of Theorem~\ref{thm:complexes1}]
  The second statement follows from the first via the bijection of
  Proposition~\ref{prop1}. To prove the first one, the remark before
  Example~\ref{ex:c23} yields an isomorphism between the face posets
  of~$\CC(r,s)$~and~$\LL^t(r,s)$. It therefore suffices to identify the latter
  complex as an interval of~$C_d^\downarrow(n)$, the poset of lower faces of
  the cyclic polytope~$C_d(n)$. This is done by the following
  inclusion-reversing maps:
  \begin{eqnarray*}
    \phi: \quad \LL^t(r,s) \to  C_d^\downarrow(n), \qquad
    \phi(\lambda) & = & \bigcap_{\mu\in M(r,s): \mu\subseteq\lambda}
    \chi^{-1}\big(\aa(\mu)\big),\\
    \psi: \quad C_d^\downarrow(n) \to 2^{[s]\times[r]}, \qquad
    \psi(G) & = & \bigcup_{F\in\FF(G)} \aa^{-1}\big(\chi(F)\big). 
  \end{eqnarray*}
  Here $\FF(G)$ denotes the set of all those  lower facets of~$C_d(n)$ that
  contain the face~$G$, and $\aa^{-1}:\{0,1,\dots,r\}^s\to M(r,s)$ the
  bijection that maps ``hole size vectors'' 
  to inclusion-minimal partial lattice paths. In particular, $\aa^{-1}\chi$ and
  $\chi^{-1}\aa$ are mutually inverse bijections between the set of lower
  facets of~$C_d(n)$ and the set $M(r,s)$.  Now note that
  \[
      \phi\psi|_{\text{Im}\phi}(G) \ = \ \phi\bigg(\bigcup_{F\in\FF(G)}
      \aa^{-1}\chi(F) \bigg) 
      \ = \ \phi\bigg(\bigcup_{\mu\in\widetilde M(G)}\mu\bigg) \ = \ 
      \bigcap_{\mu\in\widetilde M(G)} \chi^{-1}\aa(\mu) \ = \ G,
  \]
  where $\widetilde M(G)= \{\aa^{-1}\chi(F)\in M(r,s):F\supseteq G\}=
  \{\mu\in M(r,s):\chi^{-1}\aa(\mu)\supseteq G\}$, and
  that
  \[
      \psi\phi(\lambda) \ = \ \psi\bigg( \bigcap_{\mu\in M(r,s):
        \mu\subseteq\lambda} \chi^{-1}\aa(\mu) \bigg) \ = \
      \psi\bigg( \bigcap_{F:\aa^{-1}\chi(F)\subseteq\lambda} F \bigg) \ = \
      \bigcup_{F:\aa^{-1}\chi(F)\subseteq\lambda} \aa^{-1}\chi(F) \ = \ \lambda.
  \]
  It now follows that $\lambda\subset\lambda'$ in $\LL^t(r,s)$ if and only if
  $\phi(\lambda')\subset\phi(\lambda)$ in~$C_d^\downarrow(n)$. The forward
  direction is clear, as the intersection in the definition
  of~$\phi(\lambda')$ is taken over a larger subset than in~$\phi(\lambda)$,
  while the reverse direction follows from $\psi\phi(\lambda)=\lambda$ and the
  fact that $\psi$~reverses inclusions. The statement about the dimension of
  the complex follows because any member of~$M(r,s)$ has $r$~elements, while
  the cardinality of a full lattice path is~$r+s-1=r+d/2$.
\end{proof}

\subsection{Proof of Theorem~\ref{thm:subcomplexes}}

By the discussion leading up to the proof of
Theorem~\ref{thm:complexes1}, we identify the staircase
triangulation~$\Sigma(r,s)$ as $\IHom\big(K_r,S(r,s)\big)$, where
$S(r,s)$ is the graph on the vertex set $[r]\times[s]$ in which an
edge joins $(i_1,j_1)$ to $(i_2,j_2)$ exactly if $i_1<i_2$ and $j_1\le
j_2$. If we can show that $S(r,s)\subset I(k,m)$ for some values
of~$r,s$, then by the definition of~$\IHom\big(K_r,S(r,s)\big)$ and
the functoriality of $\Hom(K_r,-)$ we obtain
$\IHom\big(K_r,S(r,s)\big)\subset \Hom\big(K_r,S(r,s)\big) \subset
\Hom(K_r,H)$. In fact, $\IHom\big(K_r,S(r,s)\big) =
\Hom\big(K_r,S(r,s)\big)$ because $\omega(S(r,s))=r$.

\begin{figure}[htbp]
  \centering
  \includegraphics[width=.3\linewidth]{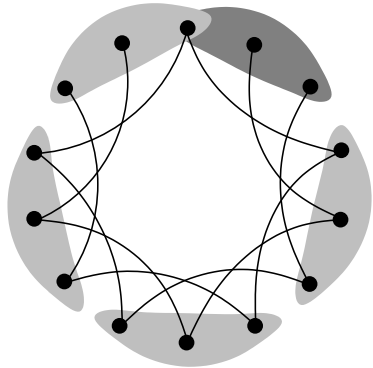}
  \caption{The construction for $k=4$, $m=6$, $r=4$ and $s=3$.}
  \label{fig:maxclique}
\end{figure}

We now find $r,s$ such that $S(r,s)\subset I(k,m)$. First, for any
$x$~and~$s$ with $0\le x\le N$ and $1\le s \le k-1$, let
$\sigma_x(s)\subset\delta(k,m)$ be the set of $s$~diagonals
\[
   \sigma_x(s) \ = \ \big\{(x+j,x+j+k-1)\bmod N:0\le j\le s-1\big\}.
\]
Clearly, any two diagonals in each $\sigma_x(s)$ cross, so that
$\sigma_x(s)$~is an independent set in~$I(k,m)$.  Now choose an
integer $r$ such that $1\le r\le \frac{(m-1)(k-2)+2}{k-1}$ and put
\[
  S \ = \ \bigcup_{b=0}^{r-1} \sigma_{b(k-1)}(s).
\]
By our choice of~$r$, the set $S$ is the vertex set of a copy
of~$S(r,s)$ inside~$I(k,m)$, as desired
(cf.~Figure~\ref{fig:maxclique}).

This proves part (a) of the theorem. Part (b) follows by combining
Theorem~\ref{prop:hom} and Proposition~\ref{prop1} with part (a). The proof of
Theorem~\ref{thm:subcomplexes} is now complete. \hfill$\qed$ 

\section{Acknowledgment}

The author would like to thank both anonymous referees for their detailed 
suggestions that helped improve the exposition.

\normalfont\unboldmath

\end{document}

%% file: graphics/cayley1.pstex_t
\begin{picture}(0,0)%
\includegraphics{cayley1.pstex}%
\end{picture}%
\setlength{\unitlength}{3947sp}%
\begingroup\makeatletter\ifx\SetFigFont\undefined%
\gdef\SetFigFont#1#2#3#4#5{%
  \reset@font\fontsize{#1}{#2pt}%
  \fontfamily{#3}\fontseries{#4}\fontshape{#5}%
  \selectfont}%
\fi\endgroup%
\begin{picture}(5408,4026)(901,-4444)
\put(5303,-3488){\makebox(0,0)[lb]{\smash{{\SetFigFont{12}{14.4}{\rmdefault}{\bfdefault}{\updefault}{\color[rgb]{0,0,0}$\pi_\Delta$}%
}}}}
\put(2251,-3488){\makebox(0,0)[lb]{\smash{{\SetFigFont{12}{14.4}{\rmdefault}{\bfdefault}{\updefault}{\color[rgb]{0,0,0}$\pi_\Delta$}%
}}}}
\put(3376,-3061){\makebox(0,0)[lb]{\smash{{\SetFigFont{12}{14.4}{\rmdefault}{\bfdefault}{\updefault}{\color[rgb]{0,0,0}$\pi_\Box\,\iota_L(\sigma)$}%
}}}}
\put(4576,-3586){\makebox(0,0)[lb]{\smash{{\SetFigFont{12}{14.4}{\rmdefault}{\bfdefault}{\updefault}{\color[rgb]{0,0,0}$L'$}%
}}}}
\end{picture}%

%% file: graphics/dimension.pstex_t
\begin{picture}(0,0)%
\includegraphics{dimension.pstex}%
\end{picture}%
\setlength{\unitlength}{3947sp}%
\begingroup\makeatletter\ifx\SetFigFont\undefined%
\gdef\SetFigFont#1#2#3#4#5{%
  \reset@font\fontsize{#1}{#2pt}%
  \fontfamily{#3}\fontseries{#4}\fontshape{#5}%
  \selectfont}%
\fi\endgroup%
\begin{picture}(5573,2576)(1166,-2667)
\put(4683,-2628){\makebox(0,0)[lb]{\smash{{\SetFigFont{11}{13.2}{\rmdefault}{\mddefault}{\updefault}{\color[rgb]{0,0,0}$x_i+\delta_i^+$}%
}}}}
\put(4092,-196){\makebox(0,0)[lb]{\smash{{\SetFigFont{11}{13.2}{\rmdefault}{\mddefault}{\updefault}{\color[rgb]{0,0,0}$y_i$}%
}}}}
\put(3171,-196){\makebox(0,0)[lb]{\smash{{\SetFigFont{11}{13.2}{\rmdefault}{\mddefault}{\updefault}{\color[rgb]{0,0,0}$y_i-\delta_i^+$}%
}}}}
\put(4552,-196){\makebox(0,0)[lb]{\smash{{\SetFigFont{11}{13.2}{\rmdefault}{\mddefault}{\updefault}{\color[rgb]{0,0,0}$y_i+\delta_i^-$}%
}}}}
\put(1791,-1050){\makebox(0,0)[lb]{\smash{{\SetFigFont{11}{13.2}{\rmdefault}{\mddefault}{\updefault}{\color[rgb]{0,0,0}$d_1$}%
}}}}
\put(5998,-985){\makebox(0,0)[lb]{\smash{{\SetFigFont{11}{13.2}{\rmdefault}{\mddefault}{\updefault}{\color[rgb]{0,0,0}$d_{m-1}$}%
}}}}
\put(4026,-2628){\makebox(0,0)[lb]{\smash{{\SetFigFont{11}{13.2}{\rmdefault}{\mddefault}{\updefault}{\color[rgb]{0,0,0}$x_i$}%
}}}}
\put(3303,-2628){\makebox(0,0)[lb]{\smash{{\SetFigFont{11}{13.2}{\rmdefault}{\mddefault}{\updefault}{\color[rgb]{0,0,0}$x_i-\delta_i^-$}%
}}}}
\put(5222,-2331){\makebox(0,0)[lb]{\smash{{\SetFigFont{11}{13.2}{\rmdefault}{\mddefault}{\updefault}{\color[rgb]{0,0,0}$x_i+\ell_i$}%
}}}}
\put(5215,-534){\makebox(0,0)[lb]{\smash{{\SetFigFont{11}{13.2}{\rmdefault}{\mddefault}{\updefault}{\color[rgb]{0,0,0}$y_i+m_i$}%
}}}}
\put(3894,-656){\makebox(0,0)[lb]{\smash{{\SetFigFont{11}{13.2}{\rmdefault}{\mddefault}{\updefault}{\color[rgb]{0,0,0}$d_i$}%
}}}}
\end{picture}%

%% file: graphics/cyclic.pstex_t
\begin{picture}(0,0)%
\includegraphics{cyclic.pstex}%
\end{picture}%
\setlength{\unitlength}{3947sp}%
\begingroup\makeatletter\ifx\SetFigFont\undefined%
\gdef\SetFigFont#1#2#3#4#5{%
  \reset@font\fontsize{#1}{#2pt}%
  \fontfamily{#3}\fontseries{#4}\fontshape{#5}%
  \selectfont}%
\fi\endgroup%
\begin{picture}(6214,2164)(2557,-2324)
\put(2745,-2282){\makebox(0,0)[lb]{\smash{{\SetFigFont{11}{13.2}{\rmdefault}{\bfdefault}{\updefault}{\color[rgb]{0,0,1}34}%
}}}}
\put(3201,-2282){\makebox(0,0)[lb]{\smash{{\SetFigFont{11}{13.2}{\rmdefault}{\bfdefault}{\updefault}{\color[rgb]{0,0,1}45}%
}}}}
\put(3656,-2282){\makebox(0,0)[lb]{\smash{{\SetFigFont{11}{13.2}{\rmdefault}{\bfdefault}{\updefault}{\color[rgb]{0,0,1}56}%
}}}}
\put(4112,-2282){\makebox(0,0)[lb]{\smash{{\SetFigFont{11}{13.2}{\rmdefault}{\bfdefault}{\updefault}{\color[rgb]{0,0,1}67}%
}}}}
\put(4510,-2282){\makebox(0,0)[lb]{\smash{{\SetFigFont{11}{13.2}{\rmdefault}{\bfdefault}{\updefault}{\color[rgb]{0,0,1}78}%
}}}}
\put(4226,-858){\makebox(0,0)[lb]{\smash{{\SetFigFont{11}{13.2}{\rmdefault}{\bfdefault}{\updefault}{\color[rgb]{1,0,0}301}%
}}}}
\put(4681,-858){\makebox(0,0)[lb]{\smash{{\SetFigFont{11}{13.2}{\rmdefault}{\bfdefault}{\updefault}{\color[rgb]{1,0,0}310}%
}}}}
\put(4681,-402){\makebox(0,0)[lb]{\smash{{\SetFigFont{11}{13.2}{\rmdefault}{\bfdefault}{\updefault}{\color[rgb]{1,0,0}400}%
}}}}
\put(4681,-1314){\makebox(0,0)[lb]{\smash{{\SetFigFont{11}{13.2}{\rmdefault}{\bfdefault}{\updefault}{\color[rgb]{1,0,0}220}%
}}}}
\put(4226,-1314){\makebox(0,0)[lb]{\smash{{\SetFigFont{11}{13.2}{\rmdefault}{\bfdefault}{\updefault}{\color[rgb]{1,0,0}211}%
}}}}
\put(3770,-1314){\makebox(0,0)[lb]{\smash{{\SetFigFont{11}{13.2}{\rmdefault}{\bfdefault}{\updefault}{\color[rgb]{1,0,0}202}%
}}}}
\put(3941,-1053){\makebox(0,0)[lb]{\smash{{\SetFigFont{11}{13.2}{\rmdefault}{\bfdefault}{\updefault}{\color[rgb]{0,0,0}46}%
}}}}
\put(5562,-1143){\makebox(0,0)[lb]{\smash{{\SetFigFont{11}{13.2}{\familydefault}{\mddefault}{\updefault}{\color[rgb]{0,0,0}2}%
}}}}
\put(5562,-2054){\makebox(0,0)[lb]{\smash{{\SetFigFont{11}{13.2}{\familydefault}{\mddefault}{\updefault}{\color[rgb]{0,0,0}1}%
}}}}
\put(5791,-2282){\makebox(0,0)[lb]{\smash{{\SetFigFont{11}{13.2}{\familydefault}{\mddefault}{\updefault}{\color[rgb]{0,0,0}1}%
}}}}
\put(6702,-2282){\makebox(0,0)[lb]{\smash{{\SetFigFont{11}{13.2}{\familydefault}{\mddefault}{\updefault}{\color[rgb]{0,0,0}2}%
}}}}
\put(7556,-2282){\makebox(0,0)[lb]{\smash{{\SetFigFont{11}{13.2}{\familydefault}{\mddefault}{\updefault}{\color[rgb]{0,0,0}3}%
}}}}
\put(8353,-2282){\makebox(0,0)[lb]{\smash{{\SetFigFont{11}{13.2}{\familydefault}{\mddefault}{\updefault}{\color[rgb]{0,0,0}$r=4$}%
}}}}
\put(5376,-301){\makebox(0,0)[lb]{\smash{{\SetFigFont{11}{13.2}{\familydefault}{\mddefault}{\updefault}{\color[rgb]{0,0,0}$s=3$}%
}}}}
\put(2557,-2101){\makebox(0,0)[lb]{\smash{{\SetFigFont{11}{13.2}{\rmdefault}{\bfdefault}{\updefault}{\color[rgb]{0,0,1}12}%
}}}}
\put(2978,-1646){\makebox(0,0)[lb]{\smash{{\SetFigFont{11}{13.2}{\rmdefault}{\bfdefault}{\updefault}{\color[rgb]{0,0,1}23}%
}}}}
\put(3426,-1194){\makebox(0,0)[lb]{\smash{{\SetFigFont{11}{13.2}{\rmdefault}{\bfdefault}{\updefault}{\color[rgb]{0,0,1}34}%
}}}}
\put(3910,-723){\makebox(0,0)[lb]{\smash{{\SetFigFont{11}{13.2}{\rmdefault}{\bfdefault}{\updefault}{\color[rgb]{0,0,1}45}%
}}}}
\put(4345,-275){\makebox(0,0)[lb]{\smash{{\SetFigFont{11}{13.2}{\rmdefault}{\bfdefault}{\updefault}{\color[rgb]{0,0,1}56}%
}}}}
\end{picture}%